\documentclass[12pt]{amsart}

\begin{document}

\title{Cohen-Macaulay modules and holonomic modules over filtered rings}

\author{Hiroki Miyahara$^{(*)}$ and Kenji Nishida$^{(**)}$}
\address{Department of Mathematical Sciences, Shinshu University,
Matsumoto, 390-8621, Japan}
\email{(*) miyahara\_shinshu@yahoo.co.jp, 
(**)kenisida@math.shinshu-u.ac.jp}
\maketitle
\begin{abstract}
We study Gorenstein dimension and grade of a module $M$ over a filtered ring whose assosiated graded ring is a commutative Noetherian ring. An equality or an inequality between these invariants of a filtered module and its associated graded module is the most valuable property for an investigation of  filtered rings. We prove an inequality G-dim$M\leq\mbox{G-dim\,gr}M$ and an equality ${\rm grade}M={\rm grade\,gr}M$, whenever Gorenstein dimension of ${\rm gr}M$ is finite (Theorems 2.3 and 2.8). We would say that the use of G-dimension adds a new viewpoint for studying filtered rings and modules. We apply these results to a filtered ring with a Cohen-Macaulay or Gorenstein associated graded ring and study a Cohen-Macaulay, perfect or holonomic module.
\end{abstract}

\keywords{}
\footnote[0]{2000 {\it Mathematics Subject Classification.} 
13C14, 13D05, 16E10, 16E30, 16E65, 16W70 \\ 
keywords: Gorenstein dimension, grade, filtered ring, Cohen-Macaulay module, holonomic module\\
(**) Research of the author is supported by Grant-in-Aid for Scientific
Researches C(2) in Japan}

\begin{center}
1. {\sc Introduction}
\end{center}

Homological theory of filtered (non-commutative) rings grew in studying , among others, $D$-modules, i.e., rings of differential operators (cf. \cite{Bj1}, \cite{K} etc.). The use of an invariant `grade' is a core of the theory for Auslander regular or Gorenstein filtered rings (\cite{Bj1}, \cite{Bj2}, \cite{Bj3}, \cite{BE}, \cite{HVO}). In particular, its invariance under forming associated graded modules is essential. Using Gorenstein dimension (\cite{AB}, \cite{C}), we extend the class of rings for which the invariance holds.

Let $\Lambda$ be a left and right Noetherian ring. Let ${\rm mod}\Lambda$ (respectively, ${\rm mod}\Lambda^{\rm op}$) be the category of all finitely generated left (respectively, right) $\Lambda$-modules. We denote the stable category by $\underline{\rm mod}\Lambda$, the syzygy functor by $\Omega:\underline{\rm mod}\Lambda \rightarrow \underline{\rm mod}\Lambda$, and the transpose functor by ${\rm Tr}:\underline{\rm mod}\Lambda \rightarrow \underline{\rm mod}\Lambda^{\rm op}$ (see \cite{ARS}, Chapter 4, \S1 or \cite{AB}, Chapter 2, \S1). For $M\in{\rm mod}\Lambda$, we put $M^*:={\rm Hom}_{\Lambda}(M,\Lambda)\in {\rm mod}\Lambda^{\rm op}$. 

Gorenstein dimension, one of the most valuable invariants of the homological study of rings and modules, is introduced in \cite{AB}. A $\Lambda$-module $M$ is said to have {\it Gorenstein dimension zero}, denoted by G-dim$_{\Lambda}M=0$, if $M^{**}\cong M$ and ${\rm Ext}_{\Lambda}^k(M,\Lambda)={\rm Ext}_{\Lambda^{\rm op}}^k(M^*,\Lambda)=0$ for $k>0$. It follows from \cite{AB}, Proposition 3.8 that G-dim $M=0$ if and only if ${\rm Ext}_{\Lambda}^k(M,\Lambda)={\rm Ext}_{\Lambda^{\rm op}}^k({\rm Tr}M,\Lambda)=0$ for $k>0$. For a positive integer $k$, $M$ is said to have {\it Gorenstein dimension less than or equal to $k$}, denoted by G-dim $M\leq k$, if there exists an exact sequence $0\rightarrow G_k\rightarrow\cdots\rightarrow G_0\rightarrow M\rightarrow 0$ with G-dim $G_i=0$ for $0\leq i\leq k$. We have that G-dim $M\leq k$ if and only if G-dim $\Omega^kM=0$ by \cite{AB}, Theorem 3.13. It is also proved in \cite{AB} that if G-dim $M<\infty$ then G-dim $M={\rm sup}\{k:{\rm Ext}_{\Lambda}^k(M,\Lambda)\neq0\}$. In the following, we abbreviate `Gorenstein dimension' to G-dimension. 

We define another important invariant `grade'. Let $M\in{\rm mod}\Lambda$. We put grade$_{\Lambda}M:={\rm inf}\{k:{\rm Ext}_{\Lambda}^k(M,\Lambda)\neq0\}$.

In this paper we study G-dimension and grade of a filtered module over a filtered ring whose assiciated graded ring is commutative and Noetherian and apply the results to a filtered ring with a Gorenstein or Cohen-Macaulay associated graded ring. 

In section two, we study G-dimension and grade of modules over a filtered ring. As usual, we analyze them by using the properties of assiciated graded modules. 
\noindent We start from studying G-dimension. When an associated graded ring gr$\Lambda$ of a filtered ring $\Lambda$ is commutative and Noetherian, a filtered $\Lambda$-module $M$ whose associated graded module gr$M$ has finite G-dimension has also finite G-dimension and an inequality G-dim$M\leq\mbox{G-dim\,gr}M$ holds true (Theorem 2.3). We see that if an associated graded ring is regular then an equality holds for every $M$. However, it is open whether an equality holds or not in general. As for G-dimension zero, we show that if G-dim\,gr$M=0$, then G-dim$M=0$ and the converse holds whenever some additional conditions for M are assumed (Theorem 2.5). Assume further that gr$\Lambda$ is a $^*$local ring with the condition (P) (see Appendix), then `Auslander-Bridger formula' holds for a filtered module $M$ such that gr$M$ has finite G-dimension and G-dim$M=\mbox{G-dim\,gr}M$:
$\mbox{G-dim}M+\,^*\! {\rm depth\,gr}M=\,^*\!{\rm depth\,gr}\Lambda$ (Proposition 2.6). 

\noindent  To handle grade in the literatures, a kind of `finitary' condition over a ring such as `regularity' or `Gorensteiness' is setted (\cite{BE}, \S5 and  \cite{HVO}, Chapter I\!I\!I, \S2, 2.5). We find out that only the finiteness of G-dimension of ${\rm gr}M$ implies ${\rm grade}M={\rm grade\,gr}M$ for a filtered module with a good filtration (Theorem 2.8). Suppose that gr$\Lambda$ is Gorenstein. Then all finite gr$\Lambda$-modules have finite G-dimension. Thus all filtered modules with a good filtration satisfy the equality. Since regularity implies Gorensteiness, our results also cover regular filtered rings.

In section three, we apply the results obtained in the previous section to Cohen-Macaulay modules over filtered rings with a Cohen-Macaulay associated graded ring and holonomic modules over Gorenstein filtered rings. When gr$\Lambda$ is a Cohen-Macaulay $^*$local ring with the condition (P), we define Cohen-Macaulay filtered modules and see that they are perfect. Then they satisfy a duality (Theorem 3.2). Moreover, assume that $\Lambda$ is Gorenstein. Then injective dimension of $\Lambda$ is finite, say $d$, so that we can define a holonomic module. A filtered module $M$ with a good filtration is holonomic, if grade$M=d$. We generalize some results in \cite{HVO}, Chapter I\!I\!I, \S4 and give a characterization of a holonomic module $M$ by a property of ${\rm Min(gr}M)$. An example of a filtered (non-regular) Gorenstein ring is given in 3.8.

The summary of commutative graded Noetherian rings, especially, $^*$local rings are stated in Appendix.

\vspace{1cm}

\begin{center}
2. {\sc Gorenstein dimension and grade for modules over filtered Noetherian rings}
\end{center}

Let $\Lambda$ be a ring. A family ${\mathcal F}=\{{\mathcal F}_p\Lambda:p\in{\Bbb N}\}$ of additive subgroups of $\Lambda$ is called a {\it filtration} of $\Lambda$, if

\medskip

\medskip

(i)    $1\in {\mathcal F}_0\Lambda,$

\medskip

(ii)   ${\mathcal F}_p\Lambda\subset{\mathcal F}_{p+1}\Lambda,$

\medskip

(iii)  $({\mathcal F}_p\Lambda)({\mathcal F}_q\Lambda)\subset{\mathcal F}_{p+q}\Lambda,$

\medskip

(iv)   $\Lambda=\cup_{p\in{\Bbb N}}{\mathcal F}_p\Lambda.$

\medskip

A pair $(\Lambda,{\mathcal F})$ is called a {\it filtered ring}. In the following, a ring $\Lambda$ is always a filtered ring for some filtration ${\mathcal F}$, so that we only say that $\Lambda$ is a filtered ring. 

Let $\sigma_p:{\mathcal F}_p\Lambda\rightarrow {\mathcal F}_p\Lambda/{\mathcal F}_{p-1}\Lambda$ be a natural homomorphism. Put
$${\rm gr}\Lambda={\rm gr}_{\mathcal F}\Lambda:=\bigoplus_{p=0}^{\infty}{\mathcal F}_p\Lambda/{\mathcal F}_{p-1}\Lambda\ \ ({\mathcal F}_{-1}\Lambda=0).$$
Then ${\rm gr}\Lambda$ is a graded ring with multiplication 
$$\sigma_p(a)\sigma_q(b)=\sigma_{p+q}(ab),\ \ a\in{\mathcal F}_p\Lambda,\ b\in{\mathcal F}_q\Lambda.$$

We always assume that ${\rm gr}\Lambda$ is a commutative Noetherian ring. Therefore, $\Lambda$ is a right and left Noetherian ring. Our main objective is to study $\Lambda$ by relating G-dimension and grade of ${\rm mod}\Lambda$ and those of ${\rm mod}({\rm gr}\Lambda)$. Sometimes we assume further that ${\rm gr}\Lambda$ is a $^*$local ring with the condition (P) (see Appendix).

Let $M$ be a (left) $\Lambda$-module. A family ${\mathcal F}=\{{\mathcal F}_pM:p\in{\Bbb Z}\}$ of additive subgroups of $M$ is called a {\it filtration} of $M$, if

\medskip

(i)    ${\mathcal F}_pM\subset{\mathcal F}_{p+1}M,$

\medskip

(ii)   ${\mathcal F}_{-p}M=0\mbox{ for }p>\!>0,$

\medskip

(iii)  $({\mathcal F}_p\Lambda)({\mathcal F}_qM)\subset{\mathcal F}_{p+q}M,$

\medskip

(iv)   $M=\cup_{p\in{\Bbb Z}}{\mathcal F}_pM.$

\medskip

A pair $(M,{\mathcal F})$ is called a {\it filtered} $\Lambda$-module. Similar to $\Lambda$, we sometimes abbreviate and say that $M$ is a filtered module. Let $\tau_p:{\mathcal F}_pM\rightarrow{\mathcal F}_pM/{\mathcal F}_{p-1}M$ be a natural homomorphism. Put
$${\rm gr}M={\rm gr}_{\mathcal F}M:=\bigoplus_{p\in{\Bbb Z}}{\mathcal F}_pM/{\mathcal F}_{p-1}M.$$

Then ${\rm gr}M$ is a graded ${\rm gr}\Lambda$-module by 
$$\sigma_p(a)\tau_q(x)=\tau_{p+q}(ax),\ \ a\in{\mathcal F}_p\Lambda,\ x\in{\mathcal F}_qM.$$
As for filtered rings and module, the reader is referred to \cite{HVO} or \cite{NVO}. We only state here some definitions and facts. For a filtered module $(M,{\mathcal F})$, we call ${\mathcal F}$ to be a {\it good filtration}, if there exist $p_k\in{\Bbb Z}$ and $m_k\in M\ (1\leq k\leq r)$ such that 
$${\mathcal F}_pM=\sum_{k=1}^r({\mathcal F}_{p-p_k}\Lambda)m_k$$
for all $p\in{\Bbb Z}$. Then the following three conditions are equivalent (\cite{HVO}, Chapter I, 5.2 and \cite{NVO}, Chapter D, IV.3)

\medskip

(a) $M$ has a good filtration.

\medskip

(b) ${\rm gr}_{\mathcal F}M$ is a finite ${\rm gr}\Lambda$-module for a filtration ${\mathcal F}$.

\medskip

(c) $M$ is a finitely generated $\Lambda$-module.

\medskip

Therefore, we only consider a good filtration for a finitely generated $\Lambda$-module $M$, so that ${\rm gr}M$ is a finite ${\rm gr}\Lambda$-module.

Let $M,\ N$ be filtered $\Lambda$-modules. A $\Lambda$-homomorphism $f:M\rightarrow N$ is called a {\it filtered homomorphism}, if $f({\mathcal F}_pM)\subset{\mathcal F}_pN$ for all $p\in{\Bbb Z}$. Further, $f$ is called {\it strict}, if $f({\mathcal F}_pM)={\rm Im}f\cap{\mathcal F}_pN$ for all $p\in{\Bbb Z}$. 
If $M'$ is a submodule of $M$, then $\{M'\cap{\mathcal F}_pM:p\in{\Bbb Z}\}$, respectively $\{{\mathcal F}_pM+M'/M':p\in{\Bbb Z}\}$ is a good filtration on $M'$, respectively $M/M'$. We call them induced filtration on $M'$ or $M/M'$ and note that the canonical homomorphisms $M'\hookrightarrow M$ and $M\rightarrow M/M'$ are strict.  

For a filtered homomorphism $f:M\rightarrow N$, we define a map $f_p:{\mathcal F}_pM/{\mathcal F}_{p-1}M\rightarrow{\mathcal F}_pN/{\mathcal F}_{p-1}N$ by $f_p(\tau_p(x))=\tau_p(f(x))$ for $x\in{\mathcal F}_pM$. Then we define a ${\rm gr}\Lambda$-homomorphism
$${\rm gr}f:{\rm gr}M=\oplus{\mathcal F}_pM/{\mathcal F}_{p-1}M\longrightarrow{\rm gr}N=\oplus{\mathcal F}_pN/{\mathcal F}_{p-1}N$$
by ${\rm gr}f:=\oplus f_p$, so that ${\rm gr}f(\tau_p(x))=\tau_p(f(x))$ for $x\in{\mathcal F}_pM$. It is easily seen that ${\rm gr}fg=({\rm gr}f)({\rm gr}g)$ for filtered homomorphisms $f:M\rightarrow N$ and $g:K\rightarrow M$. 

For a filtered module $M$, an exact sequence 
$$\cdots \longrightarrow F_i \overset{f_i}{\longrightarrow}\cdots \overset{f_1}{\longrightarrow}F_0\overset{f_0}{\longrightarrow}M\longrightarrow0$$
is called a {\it filtered free resolution} of $M$, if all $F_i$ are filtered free $\Lambda$-modules and all homomorphisms are strict filtered homomorphisms. We can always constract such a resolution with all $F_i$ of finite rank for a finitely generated $\Lambda$-module (see \cite{NVO}, Chapter D, IV). 

Let $M,\ N$ be filtered $\Lambda$-modules. We put, for $p\in{\Bbb Z}$,
$${\mathcal F}_p{\rm Hom}_{\Lambda}(M,N)=\{f\in{\rm Hom}_{\Lambda}(M,N):f({\mathcal F}_qM)\subset{\mathcal F}_{p+q}N\mbox{ for all }q\in{\Bbb Z}\}$$
Then we have an ascending chain
$$\cdots\subset{\mathcal F}_p{\rm Hom}_{\Lambda}(M,N)\subset{\mathcal F}_{p+1}{\rm Hom}_{\Lambda}(M,N)\subset\cdots$$
of additive subgroups of ${\rm Hom}_{\Lambda}(M,N)$. Set
$${\rm gr\,Hom}_{\Lambda}(M,N):=\bigoplus_{p\in{\Bbb Z}}{\mathcal F}_p{\rm Hom}_{\Lambda}(M,N)/{\mathcal F}_{p-1}{\rm Hom}_{\Lambda}(M,N)$$
Define an additive homomorphism
$$\varphi=\varphi(M,N):{\rm gr\,Hom}_{\Lambda}(M,N)\longrightarrow{\rm Hom}_{{\rm gr}\Lambda}({\rm gr}M,{\rm gr}N),\ \varphi(\tau_p(f))(\tau_q(x))=\tau_{p+q}(f(x))$$
for $f\in{\mathcal F}_p{\rm Hom}_{\Lambda}(M,N),\ x\in{\mathcal F}_qM$, where 
$$\tau_p:{\mathcal F}_p{\rm Hom}_{\Lambda}(M,N)\longrightarrow{\mathcal F}_p{\rm Hom}_{\Lambda}(M,N)/{\mathcal F}_{p-1}{\rm Hom}_{\Lambda}(M,N)$$
is a natural homomorphism for every $p\in{\Bbb Z}$. When $M$ is a filtered module with a good filtration, the following facts hold (see \cite{HVO}, Chapter I, 6.9 or \cite{NVO}, Chapter D, VI.6):

\medskip

(1) ${\rm Hom}_{\Lambda}(M,N)=\cup_{p\in{\Bbb Z}}{\mathcal F}_p{\rm Hom}_{\Lambda}(M,N).$

\medskip

(2) ${\mathcal F}_{-p}{\rm Hom}_{\Lambda}(M,N)=0$ for $p>\!>0.$

\medskip

(3) $\varphi$ is injective. Moreover, if $M$ is a filtered free module, then it is bijective.

\medskip

(4) When $N=\Lambda$, an additive group ${\rm Hom}_{\Lambda}(M,\Lambda)$ is a filtered $\Lambda^{\rm op}$-module with a good filtration ${\mathcal F}:=\{{\mathcal F}_p{\rm Hom}_{\Lambda}(M,\Lambda):p\in{\Bbb Z}\}$ and $\varphi$ is a ${\rm gr}\Lambda$-homomorphism.

\medskip

Let $M\overset{f}{\rightarrow}N\overset{g}{\rightarrow}K$ be an exact sequence of filtered modules and filtered homomorphisms. Then ${\rm gr}M\overset{{\rm gr}f}{\rightarrow}{\rm gr}N\overset{{\rm gr}g}{\rightarrow}{\rm gr}K$ is exact (in ${\rm mod\,gr}\Lambda$) if and only if $f$ and $g$ are strict (see \cite{HVO}, Chapter I, 4.2.4 or \cite{NVO}, Chapter D, I\!I\!I.3).

\medskip

The following proposition is well-known.

\medskip

2.1. {\sc Proposition}.\ \ 
{\sl Let $M$ be a filtered $\Lambda$-module with a good filtation. Then ${\rm gr\,Ext}_{\Lambda}^i(M,\Lambda)$ is a subfactor of ${\rm Ext}_{{\rm gr}\Lambda}^i({\rm gr}M,{\rm gr}\Lambda)$ for $i\geq 0$.  }

{\it Proof.} See \cite{Bj1}, Chapter 2, 6.10 or \cite{HVO}, Chapter I\!I\!I, 2.2.4. $\square$

\medskip

When G-dim${\, \rm gr}M=0$, the functor Tr commutes with associated gradation. 

\medskip

2.2. {\sc Lemma}.\ \ 
{\sl Let $M$ be a filtered $\Lambda$-module with a good filtation. Then there exists an epimorphism $\alpha:{\rm Tr}_{{\rm gr}\Lambda}({\rm gr}M)\rightarrow{\rm gr}({\rm Tr}_{\Lambda}M)$.

Moreover, if {\rm G}-${\rm dim\,gr}M=0$ or {\rm grade\,gr}$M>1$, then $\alpha$ is an isomorphism.}

{\it Proof.} Take a filtered free resolution of $M$:
$$\cdots \longrightarrow F_2 \overset{f_2}{\longrightarrow}F_1\overset{f_1}{\longrightarrow}F_0\overset{f_0}{\longrightarrow}M\longrightarrow0.$$
By definition, we have an exact sequence
$$F_0^* \overset{f_1^*}{\longrightarrow}F_1^*\overset{g}{\longrightarrow}{\rm Tr}_{\Lambda}M={\rm Cok}f_1^* \longrightarrow0,$$
where $g$ is a canonical epimorphism. Let ${\rm Tr}_{\Lambda}M$ be equipped with the induced filtration by $g$. Then $g$ is a strict filtered epimorphism. Let us consider the following diagrams in ${\rm mod\,gr}\Lambda$ with the commutative squares and all the $\varphi$'s isomorphisms:

$$\begin{array}{cccccccl}
&{\rm Hom}_{{\rm gr}\Lambda}({\rm gr}F_0,{\rm gr}\Lambda)&\overset{({\rm gr}f_1)^*}{\longrightarrow}&{\rm Hom}_{{\rm gr}\Lambda}({\rm gr}F_1,{\rm gr}\Lambda)&\longrightarrow&{\rm Tr}_{{\rm gr}\Lambda}({\rm gr}M)&\longrightarrow&0(\mbox{exact})\\
(1)&\varphi\uparrow&&\varphi\uparrow&&&&\\
&{\rm gr}F_0^*&\overset{{\rm gr}(f_1^*)}{\longrightarrow}&{\rm gr}F_1^*&\overset{{\rm gr}g}{\longrightarrow}&{\rm gr}({\rm Tr}_{\Lambda}M)&\longrightarrow&0
\end{array} $$

$$\begin{array}{cccccc}
&{\rm Hom}_{{\rm gr}\Lambda}({\rm gr}F_0,{\rm gr}\Lambda)&\overset{({\rm gr}f_1)^*}{\longrightarrow}&{\rm Hom}_{{\rm gr}\Lambda}({\rm gr}F_1,{\rm gr}\Lambda)&\overset{({\rm gr}f_2)^*}{\longrightarrow}&{\rm Hom}_{{\rm gr}\Lambda}({\rm gr}F_2,{\rm gr}\Lambda)\\
(2)&\varphi\uparrow&&\varphi\uparrow&&\varphi\uparrow\\
&{\rm gr}F_0^*&\overset{{\rm gr}(f_1^*)}{\longrightarrow}&{\rm gr}F_1^*&\overset{{\rm gr}(f_2^*)}{\longrightarrow}&{\rm gr}F_2^*
\end{array} $$
Since the induced sequence $\cdots \rightarrow{\rm gr}F_1\overset{{\rm gr}f_1}{\longrightarrow}{\rm gr}F_0\rightarrow {\rm gr}M\rightarrow 0$ is a free resolution of ${\rm gr}M$, the first row of (1) is exact. Since $g$ is strict, ${\rm gr}g$ is surjective. Hence there exists a graded epimorphism $\alpha:{\rm Tr}_{{\rm gr}\Lambda}({\rm gr}M)\rightarrow{\rm gr}({\rm Tr}_{\Lambda}M)$. By assumption, we see that ${\rm Ext}^1_{{\rm gr}\Lambda}({\rm gr}M,{\rm gr}\Lambda)=0$, so that the first row of (2) is exact. There exists a filtered homomorphism $h:{\rm Tr}_{\Lambda}M\rightarrow F_2^*$ such that $f_2^*=h\circ g$. Since ${\rm gr}f_2^*={\rm gr}h\circ{\rm gr}g$, we have ${\rm Im\,gr}f_1^*\subset{\rm Ker\,gr}g\subset{\rm Ker\,gr}f_2^*$. The exactness of the second row of (2) implies ${\rm Im\,gr}f_1^*={\rm Ker\,gr}f_2^*$. Thus ${\rm Im\,gr}f_1^*={\rm Ker\,gr}g$, hence the second row of (1) is also exact, which implies that $\alpha$ is an isomorphism. $\square$

\medskip

2.3. {\sc Theorem.}\ \ 
{\sl Let $M$ be a filtered $\Lambda$-module with a good filtration such that ${\rm gr}M$ is of finite G-dimension. Then G-dim$M\leq \mbox{G-dim\,gr}M$.}

\medskip

{\it Proof.}
We show that if G-dim\,gr$M=k<\infty$, then G-dim$M\leq k$. Let $k=0$. Assume that G-dim${\,\rm gr}M=0$. For $i>0$, since ${\rm gr\,Ext}_{\Lambda}^i(M,\Lambda)$ is a subfactor of ${\rm Ext}_{{\rm gr}\Lambda}^i({\rm gr}M,{\rm gr}\Lambda)$, we have ${\rm gr\,Ext}_{\Lambda}^i(M,\Lambda)=0$. Hence ${\rm Ext}_{\Lambda}^i(M,\Lambda)=0$. By Lemma 2.2, ${\rm Ext}_{{\rm gr}\Lambda}^i({\rm gr\,Tr}_{\Lambda}M,{\rm gr}\Lambda)\cong  {\rm Ext}_{{\rm gr}\Lambda}^i({\rm Tr}_{{\rm gr}\Lambda}({\rm gr}M),{\rm gr}\Lambda)=0$ for $i>0$. Hence ${\rm Ext}_{\Lambda^{\rm op}}^i({\rm Tr}_{\Lambda}M,\Lambda)=0$ as above. Thus G-${\rm dim}M=0$.

Let $k>0$. Since ${\rm gr}(\Omega^kM)$ and $\Omega^k({\rm gr}M)$ are stably isomorphic(see \cite{F}, p.226 for the definition), the following holds: 
$$\mbox{G-dim}\,{\rm gr}M\leq k\Leftrightarrow \mbox{G-dim}\,\Omega^k({\rm gr}M)=0\Leftrightarrow \mbox{G-dim}\,{\rm gr}(\Omega^kM)=0.$$
Thus the statement holds by the case of $k=0$. $\square$

\medskip

2.4. {\sc Corollary.}\ \ 
{\sl Assume that ${\rm gr}\Lambda$ is a $^*$local ring with the condition (P). If ${\rm gr}\Lambda$ is Gorenstein, then ${\rm id}_{\Lambda}\Lambda={\rm id}_{\Lambda^{\rm op}}\Lambda\leq\,^*{\rm depth\,gr}\Lambda$. }

{\it Proof.}
Let $M$ be a finitely generated $\Lambda$-module. Then $M$ is a filtered module with a good filtration. Then G-dim$\,{\rm gr}M<\infty$ by Theorem A.9. Hence 
$$\mbox{G-dim}M\leq \mbox{\rm G-dim\,gr}M=\,^*{\rm depth\,gr}\Lambda-\,^*{\rm depth\,gr}M\leq\,^*{\rm depth\,gr}\Lambda.$$
Therefore, ${\rm Ext}_{\Lambda}^i(M,\Lambda)=0$ for all $i>\,^*{\rm depth\,gr}\Lambda$, so that ${\rm id}_{\Lambda}\Lambda<\infty$. Similarly, we have ${\rm id}_{\Lambda^{\rm op}}\Lambda<\infty$. Thus ${\rm id}_{\Lambda}\Lambda={\rm id}_{\Lambda^{\rm op}}\Lambda\leq\,^*{\rm depth\,gr}\Lambda.$ $\square$

\medskip

Thanks to Corollary 2.4, we call a filtered ring $\Lambda$ a "{\it Gorenstein filtered ring}", if ${\rm gr}\Lambda$ is a Gorenstein $^*$local ring with the condition (P). 

\medskip

We give a necessary and sufficient condition when G-dim$\,{\rm gr}M=0$.

\medskip

2.5. {\sc Theorem.}\ \ 
{\sl Let $M$ be a filtered $\Lambda$-module with a good filtration. Then the following (1) and (2) are equivalent.

\medskip

(1) {\rm G-dim}$\,{\rm gr}M=0$.

\medskip

(2) (2.1) {\rm G-dim}$M=0$.

\medskip

\ \ \ \ \ (2.2) Suppose that $\cdots\rightarrow F_1\overset{f_1}{\longrightarrow}F_0\overset{f_0}{\longrightarrow}M\rightarrow0$ is a filtered free resolution of $M$, then all $f_i^*\ (i>0)$ are strict.

\medskip

\ \ \ \ \ (2.2\,$^*$) Suppose that $\cdots\rightarrow G_1\overset{g_1}{\longrightarrow}G_0\overset{g_0}{\longrightarrow}M^*\rightarrow0$ is a filtered free resolution of $M^*$, then all $g_i^*\ (i>0)$ are strict.

\medskip

\ \ \ \ \ (2.3) A canonical map $\theta:M\rightarrow M^{**}$ is strict. 

\medskip

Moreover, under the above conditions, 
$\varphi_M:{\rm gr}M^*\rightarrow({\rm gr}M)^*$ and $\varphi_{M^*}:{\rm gr}M^{**}\rightarrow({\rm gr}M^*)^*$ are isomorphisms, where $\varphi_M=\varphi(M,\Lambda),\ \varphi_{M^*}=\varphi(M^*,\Lambda)$.   }

\medskip

{\it Proof.}
(1) $\Rightarrow$ (2): It follows from Theorem 2.3 that G-dim$M=0$. From a filtered free resolution of $M$ in (2.2), we get an exact sequence
$$0\longrightarrow M^*\overset{f_0^*}{\longrightarrow}F_0^*\overset{f_1^*}{\longrightarrow}F_1^*\longrightarrow\cdots$$
This exact sequence and an exact sequence in ${\rm mod\,gr}\Lambda$:
$$\cdots\longrightarrow{\rm gr}F_1\longrightarrow {\rm gr}F_0\longrightarrow{\rm gr}M\longrightarrow0$$
induced from a resolution in (2.2) give the following commutative diagram
$$\begin{array}{cccccccccc}
&0&\longrightarrow&{\rm gr}M^*&\overset{{\rm gr}(f_0^*)}{\longrightarrow}&{\rm gr}F_0^*&\overset{{\rm gr}(f_1^*)}{\longrightarrow}&{\rm gr}F_1^*&\longrightarrow&\cdots \\
(*)&&&\varphi\downarrow&&\varphi_0\downarrow&&\varphi_1\downarrow&&\\
&0&\longrightarrow&({\rm gr}M)^*&\longrightarrow&({\rm gr}F_0)^*&\longrightarrow&({\rm gr}F_1)^*&\longrightarrow&\cdots
\end{array}$$
where $\varphi=\varphi(M,\Lambda),\ \varphi_i=\varphi(F_i,\Lambda)$. Since G-dim$\,{\rm gr}M=0$, the second row is exact. For $i\geq 0$, $\varphi_i$ are isomorphisms. Thus a sequence
$${\rm gr}F_0^*\overset{{\rm gr}(f_1^*)}{\longrightarrow}{\rm gr}F_1^*\overset{{\rm gr}(f_2^*)}{\longrightarrow}{\rm gr}F_2^*\longrightarrow\cdots$$
is exact, and so $f_1^*,\ f_2^*,\cdots$ are strict. Hence (2.2) holds. Since $f_0$ is a strict filtered epimorphism, $f_0^*$ is a strict filtered monomorphism. Thus the first row of $(*)$ is exact. Therefore, $\varphi:{\rm gr}M^*\rightarrow({\rm gr}M)^*$ is an isomorphism. Since G-dim$({\rm gr}M)^*=0$, we have G-dim$\,{\rm gr}M^*=0$. Hence (2.2$^*$) holds and $\varphi_{M^*}$ is an isomorphism.

Let $\eta:{\rm gr}M\rightarrow({\rm gr}M)^{**}$ be a canonical homomorphism. Consider the commutative diagram
$$\begin{array}{cccc}
&{\rm gr}M&\overset{{\rm gr}\theta}{\longrightarrow}&{\rm gr}M^{**}\\
(**)\ \ \ &\eta\downarrow&&\downarrow\,\varphi_{M^*}\\
&({\rm gr}M)^{**}&\overset{\varphi_M^*}{\longrightarrow}&({\rm gr}M^*)^*.
\end{array}$$
Since $\eta,\ \varphi_M^*,\ \varphi_{M^*}$ are isomorphisms, ${\rm gr}\theta$ is also an isomorphism. Thus $\theta$ is strict.

(2) $\Rightarrow$ (1): By (2.1) and (2.2), the first row of the diagram $(*)$ is exact. Thus the second row of $(*)$ is exact, so that ${\rm Ext}_{{\rm gr}\Lambda}^i({\rm gr}M,{\rm gr}\Lambda)=0$ for $i>0$ and $({\rm gr}M)^*\cong {\rm gr}M^*$. Since G-dim$M^*=0$, using the diagram $(*)$ obtained from (2.2$^*$), we can show that ${\rm Ext}_{{\rm gr}\Lambda}^i({\rm gr}M^*,{\rm gr}\Lambda)=0$ for $i>0$ and $({\rm gr}M^*)^*\cong{\rm gr}M^{**}$. Thus we have ${\rm Ext}_{{\rm gr}\Lambda}^i(({\rm gr}M)^*,{\rm gr}\Lambda)=0$ for $i>0$. By (2.3) and the above argument, the maps ${\rm gr}\theta,\ \varphi_M^*$ and $\varphi_{M^*}$ are isomorphisms in the diagram $(**)$, so that $\eta$ is an isomorphism. Thus G-dim$\,{\rm gr}M=0.$ $\square$
 
 \medskip
 
Let fil$\Lambda$ be a category of all filtered $\Lambda$-modules with a good filtration and filtered homomorphisms. Let $\mathcal G$ be a subcategory of fil$\Lambda$ consisting of all filtered modules $M$ whose associated graded module gr$M$ has finite G-dimension. It holds from Theorem 2.3 that a module in $\mathcal G$ has finite G-dimension. We further put a subcategory ${\mathcal G}_e$ of $\mathcal G$
$${\mathcal G}_e:=\{M\in{\mathcal G}:\mbox{G-dim}M=\mbox{G-dim\,gr}M \mbox{ for some good filtration of } M\}.$$

\medskip 
 
2.6. {\sc Proposition}.\ \ 
{\sl Assume that ${\rm gr}\Lambda$ is a $^*$local ring with the condition (P).
  Let $M\in{\mathcal G}_e$. Then the following equality holds.
$$\mbox{\rm G-dim}M+\,^*{\rm depth\,gr}M=\,^*{\rm depth\,gr}\Lambda.$$}

{\it Proof.}
The statement follows from Theorem A.8. $\square$

\medskip

2.7. {\sc Remarks}.\ \ 
(i) It is interesting to know when ${\mathcal G}_e={\mathcal G}$. If this is true, then we see that G-dim$M=0$ if and only if G-dim\,gr$M=0$ for $M\in{\mathcal G}$. Hence the condition (2.2), (2.2$^*$), (2.3) in Theorem 2.5 are superfluous.
\medskip

(ii) Suppose that $0\rightarrow M'\rightarrow M\rightarrow M''\rightarrow 0$ is a strict exact sequence of fil$\Lambda$. Then the followings are easy consequence of \cite{C}, Corollary 1.2.9 (b).

If $M', M''\in {\mathcal G}_e$ and G-dim$M'>\mbox{G-dim}M''$, then $M\in{\mathcal G}_e$.

If $M, M''\in {\mathcal G}_e$ and G-dim$M>\mbox{G-dim}M''$, then $M'\in{\mathcal G}_e$.

\medskip

 We shall study the another valuable invariant `grade'. Its nicest feature that an equation ${\rm grade}_{\Lambda}M={\rm grade}_{{\rm gr}\Lambda}{\rm gr}M$ holds for a good filtered $\Lambda$-module $M$ is proved when ${\rm gr}\Lambda$ is regular (see e.g. \cite{HVO}). We prove this equation under `module-wise' conditions by which we can apply this equation fairly wide classes of filtered rings.

\medskip
 
2.8. {\sc Theorem.}\ \ 
{\sl Let $\Lambda$ be a filtered ring such that ${\rm gr}\Lambda$ is a commutative Noetherian ring and $M$ a filtered $\Lambda$-module with a good filtration. Assume that ${\rm gr}M$ has finite G-dimension. Then an equality ${\rm grade}_{\Lambda}M={\rm grade}_{{\rm gr}\Lambda}{\rm gr}M$ holds.  }

{\it Proof.}
Put $s={\rm grade}_{{\rm gr}\Lambda}{\rm gr}M$. In order to show that grade$_{\Lambda}M=s$, we must prove:

\medskip

(i)  ${\rm Ext}_{\Lambda}^s(M,\Lambda)\neq0$,

\medskip

(ii) ${\rm Ext}_{\Lambda}^i(M,\Lambda)=0$ for $i<s$.

\medskip

2.8.1. (cf. \cite{HVO}, Chapter I\!I\!I, \S1)\ \ 
Let $\cdots \rightarrow F_i\overset{f_i}{\rightarrow}\cdots\overset{f_1}{\rightarrow}F_0\overset{f_0}{\rightarrow}M\rightarrow0$ be a filtered free resolution of $M$. Applying $(-)^*$ to it, we get a complex
$$F_{\bullet}:0\rightarrow F_0^*\overset{f_1^*}{\rightarrow}\cdots \rightarrow F_{i-2}^*\overset{f_{i-1}^*}{\rightarrow}F_{i-1}^*\overset{f_i^*}{\rightarrow}F_i^* \rightarrow\cdots$$
with each $F_i^*$ filtered free and $f_i^*$ a filtered homomorphism. We put, for $p,r,i\in{\Bbb N},$
$$Z_p^r(i):=(f_i^*)^{-1}({\mathcal F}_{p-r}F_i^*)\cap{\mathcal F}_pF_{i-1}^*,\ \ \ Z_p^{\infty}(i):={\rm Ker}f_i^*\cap{\mathcal F}_pF_{i-1}^*,$$
$$B_p^r(i):=f_{i-1}^*({\mathcal F}_{p+r-1}F_{i-2}^*)\cap{\mathcal F}_pF_{i-1}^*,\ \ \ B_p^{\infty}(i):={\rm Im}f_{i-1}^*\cap{\mathcal F}_pF_{i-1}^*$$
Then the following sequence of inclusions holds:
$$Z_p^0(i)\supset Z_p^1(i)\supset\cdots\supset Z_p^{\infty}(i)\supset B_p^{\infty}(i)\supset\cdots \supset B_p^1(i)\supset B_p^0(i).$$
We put
$$E_p^r(i):=\frac{Z_p^r(i)+{\mathcal F}_{p-1}F_{i-1}^*}{B_p^r(i)+{\mathcal F}_{p-1}F_{i-1}^*},\ \ \ E_i^r:=\bigoplus_pE_p^r(i).$$   
Then $E_i^r$ is a ${\rm gr}\Lambda$-module for $r,i\geq0$. When $r=0$, we have
$$E_i^0=\bigoplus_p\frac{(f_i^*)^{-1}({\mathcal F}_pF_i^*)\cap{\mathcal F}_pF_{i-1}^*+{\mathcal F}_{p-1}F_{i-1}^*}{f_{i-1}^*({\mathcal F}_{p-1}F_{i-2}^*)\cap{\mathcal F}_pF_{i-1}^*+{\mathcal F}_{p-1}F_{i-1}^*}=\bigoplus_p\frac{{\mathcal F}_pF_{i-1}^*}{{\mathcal F}_{p-1}F_{i-1}^*}={\rm gr}F_{i-1}^*.$$
Hence we get a complex
$$E_{\bullet}^0:0\rightarrow{\rm gr}F_0^*\rightarrow\cdots\rightarrow{\rm gr}F_i^*\rightarrow\cdots$$
which is an associated graded complex of $F_{\bullet}$. We show, for $r\geq1$, that $\{E_i^r\}_{i\geq0}$ also gives a complex $E_{\bullet}^r$. To do so, we define morphisms. By computation, it holds that 
$$E_p^r(i)=\frac{Z_p^r(i)}{B_p^r(i)+Z_{p-1}^{r-1}(i)},\ \ \ f_i^*(Z_p^r(i))={\mathcal F}_{p-r}F_i^*\cap f_i^*({\mathcal F}_pF_{i-1}^*)=B_{p-r}^{r+1}(i+1).$$
Thus the following hold:

$(1)\ \ f_i^*(Z_p^r(i))=B_{p-r}^{r+1}(i+1)\subset Z_{p-r}^r(i+1),$

$(2)\ \ f_i^*(B_p^r(i))=0\mbox{ and }f_i^*(Z_{p-1}^{r-1}(i))=B_{p-r}^r(i+1).$

\medskip

\noindent We can show that $f_i^*$ induces a map
$\tilde{f}_p^r(i):E_p^r(i)\rightarrow E_{p-r}^r(i+1),$ by

$\tilde{f}_p^r(i)(x+B_p^r(i)+{\mathcal F}_{p-1}F_{i-1}^*)=f_i^*(x)+B_{p-r}^r(i+1)+Z_{p-r-1}^{r-1}(i+1)\ (x\in Z_p^r(i)).$

\noindent Hence $\tilde{f}_p^r(i)(p\in{\Bbb N})$ give a graded ${\rm gr}\Lambda$-homomorphism
$$\tilde{f}_i^r:E_i^r=\bigoplus_pE_p^r(i)\longrightarrow E_{i+1}^r=\bigoplus_pE_p^r(i+1)$$
of degree $-r$. It is easily seen that $E_{\bullet}^r:\cdots\rightarrow E_i^r\overset{\tilde{f}_i^r}{\longrightarrow}E_{i+1}^r\rightarrow\cdots$ is a complex.

\medskip

2.8.2. {\sc Lemma}. (cf. \cite{HVO}, p.130 (6))\ \ 
{\sl Under the above notation, we have $H^i(E_{\bullet}^r)\cong E_i^{r+1}$.}

{\it Proof.} 
We show 
$$H(E_{p+r}^r(i-1)\overset{f}{\longrightarrow}E_p^r(i)\overset{g}{\longrightarrow}E_{p-r}^r(i+1))\cong E_p^{r+1}(i),$$
where we put $f:=\tilde{f}_{p+r}^r(i-1),\ g:=\tilde{f}_p^r(i)$. Using (1) and (2), we can show that
$$x+B_p^r(i)+{\mathcal F}_{p-1}F_{i-1}^*\in{\rm Ker}g\Longleftrightarrow x\in (f_i^*)^{-1}(B_{p-r}^r(i+1)+{\mathcal F}_{p-r-1}F_i^*).$$
Thus we get
$${\rm Ker}g=\frac{(Z_p^r(i)\cap(f_i^*)^{-1}(B_{p-r}^r(i+1)+{\mathcal F}_{p-r-1}F_i^*)+{\mathcal F}_{p-1}F_{i-1}^*}{B_p^r(i)+{\mathcal F}_{p-1}F_{i-1}^*}.$$
Further, we have
$${\rm Im}f=\frac{f_{i-1}^*(Z_{p+r}^r(i-1))+{\mathcal F}_{p-1}F_{i-1}^*}{B_p^r(i)+{\mathcal F}_{p-1}F_{i-1}^*}.$$
Hence the desired homology is
{\allowdisplaybreaks
\begin{align*}
\frac{{\rm Ker}g}{{\rm Im}f}
& = \frac{(Z_p^r(i)\cap(f_i^*)^{-1}(B_{p-r}^r(i+1)+{\mathcal F}_{p-r-1}F_i^*)+{\mathcal F}_{p-1}F_{i-1}^*}{f_{i-1}^*(Z_{p+r}^r(i-1))+{\mathcal F}_{p-1}F_{i-1}^*}\\
& = \frac{Z_{p-1}^{r-1}(i)+Z_p^r(i)\cap(f_i^*)^{-1}({\mathcal F}_{p-r-1}F_i^*)+{\mathcal F}_{p-1}F_{i-1}^*}{f_{i-1}^*(Z_{p+r}^r(i-1))+{\mathcal F}_{p-1}F_{i-1}^*}\\
& = \frac{Z_{p-1}^{r-1}(i)+Z_p^{r+1}(i)+{\mathcal F}_{p-1}F_{i-1}^*}{f_{i-1}^*(Z_{p+r}^r(i-1))+{\mathcal F}_{p-1}F_{i-1}^*}\\
& = \frac{Z_p^{r+1}(i)+{\mathcal F}_{p-1}F_{i-1}^*}{B_p^{r+1}(i)+{\mathcal F}_{p-1}F_{i-1}^*}\\
& = E_p^{r+1}(i),
\end{align*}}
where (2) (respectively, (1)) is used to show the second (respectively, fourth) equality.
 $\square$

\medskip

2.8.3. {\sc Corollary}.\ \ 
{\sl Assume that $E_{i-1}^1=0$. Then we have $E_{i-1}^r=0$ for $r\geq 1$ and there exists an exact sequence
$$0\rightarrow E_i^{r+1}\rightarrow E_i^r\rightarrow E_{i+1}^r$$
of ${\rm gr}\Lambda$-modules for each $r\geq1$.}

{\it Proof.} 
The first assertion directly follows from Lemma 2.8.2. Then the complex $E_{\bullet}^r$ yields an exact sequence
$0\rightarrow H^i(E_{\bullet}^r)\rightarrow E_i^r\rightarrow E_{i+1}^r.$ Since $H^i(E_{\bullet}^r)\cong E_i^{r+1}$ by lemma 2.8.2, we get the desired exact sequence. $\square$

\medskip

2.8.4.\ \ 
We will show in this subsection that $E_{s+1}^r\neq 0$.

Condider the following commutative diagram
$$\begin{array}{cccccccccc}
E_{\bullet}^0={\rm gr}(F_{\bullet}^{\ast}): & 0&\rightarrow&{\rm gr}F_0^*&\rightarrow&\cdots&\rightarrow&{\rm gr}F_i^*&\rightarrow&\cdots \\
&&&\wr|&&&&\wr|&&\\
&0&\rightarrow&({\rm gr}F_0)^*&\rightarrow&\cdots&\rightarrow&({\rm gr}F_i)^*&\rightarrow&\cdots,
\end{array}$$
where rows are complexes and the second row is obtained by applying ${\rm Hom}_{\rm gr \Lambda}(-,{\rm gr}\Lambda)$ to a free resolution
$\cdots\rightarrow{\rm gr}F_1\rightarrow{\rm gr}F_0\rightarrow{\rm gr}M\rightarrow0$ of ${\rm gr}M$. Hence an isomorphism $E_{i+1}^1\cong{\rm Ext}_{\rm gr\Lambda}^i({\rm gr}M,{\rm gr}\Lambda)$ holds by Lemma 2.8.2. (Note that $E_i^0\cong{\rm gr}F_{i-1}^*$.)

By assumption, we can apply A.15 to ${\rm gr}M$ and get the fact that ${\rm grade\,Ext}_{{\rm gr}\Lambda}^s({\rm gr}M,{\rm gr}\Lambda)=s$.  Hence it holds that grade$E_{s+1}^1=s$ and $E_{i+1}^1=0$ for $i<s$. By Corollary 2.8.3, we get an exact sequence of ${\rm gr}\Lambda$-modules
$$(3)\ \ \ 0\rightarrow E_{s+1}^{r+1}\rightarrow E_{s+1}^r\overset{\varphi}{\longrightarrow}E_{s+2}^r.$$

\noindent By Lemma 2.8.2, $E_{s+2}^r$ is a subfactor of $E_{s+2}^{r-1}$ for $r\geq 1$. Thus every ${\rm gr}\Lambda$-submodule $U$ of $E_{s+2}^r$ is also a subfactor of $E_{s+2}^1={\rm Ext}_{{\rm gr}\Lambda}^{s+1}({\rm gr}M,{\rm gr}\Lambda$), so that there exist ${\rm gr}\Lambda$-submodules $X,Y\subset{\rm Ext}_{{\rm gr}\Lambda}^{s+1}({\rm gr}M,{\rm gr}\Lambda$) such that $U\cong X/Y$. Since ${\rm grade}X\geq s+1$ and ${\rm grade}Y\geq s+1$ by A.14, it holds that ${\rm grade}U\geq s+1$. Therefore, ${\rm grade}({\rm Im}\varphi_r)\geq s+1$ for $r\geq 1$. Consider the exact sequence induced from (3):

$$0\rightarrow E_{s+1}^{r+1}\rightarrow E_{s+1}^r\rightarrow {\rm Im}\varphi_r\rightarrow0.$$
Assume that ${\rm grade}E_{s+1}^r=s$. Then ${\rm grade}E_{s+1}^{r+1}=s$ holds. Hence ${\rm grade}E_{s+1}^r=s$ holds for all $r\geq 1$ by induction. Especially, $E_{s+1}^r\neq0$ holds for all $r\geq1$.

\medskip

2.8.5. {\sc Lemma}\ \ 
{\sl There is an isomorphism $E_{i+1}^r\cong {\rm gr}({\rm Ext}_{\Lambda}^i(M,\Lambda))$ for $i\geq0$ and $r>\!>0$.}

{\it Proof}. 
Since the filtration ${\mathcal F}$ of $\Lambda$ is Zariskian (see \cite{HVO}, Chapter I, \S2, 2.4; \S3, 3.3 and Chapter I\!I, \S2, 2.1, and Proposition 2.2.1), the lemma follows from \cite{HVO}, Chapter I\!I\!I, \S2, Lemma 2.2.1(p. 150) and \S1, Corollary 1.1.7(p. 133). $\square$

\medskip

2.8.6. \ \ 
We have shown that $E_{s+1}^r\neq0$. Hence ${\rm Ext}_{\Lambda}^s(M,\Lambda)\neq0$ by Lemma 2.8.5. Therefore, (i) holds. 

Conversely, since grade$\,{\rm gr}M=s$, we have ${\rm Ext}_{\rm gr \Lambda}^i({\rm gr}M,{\rm gr}\Lambda)=0$ for $i<s$. Since ${\rm gr\,Ext}_{\Lambda}^i(M,\Lambda)$ is a subfactor of ${\rm Ext}_{\rm gr\Lambda}^i({\rm gr}M,{\rm gr \Lambda})$ by Proposition 2.1, we have ${\rm gr\,Ext}_{\Lambda}^i(M,\Lambda)=0$ for $i<s$. Therefore, ${\rm Ext}_{\Lambda}^i(M,\Lambda)=0$ for $i<s$, so that (ii) holds. This accomplishes the proof of 2.8. $\square$

\medskip

2.9. {\sc Remarks.}\ \ 
(i) Let $M\in{\mathcal G}$. Then it follows from 2.3 and 2.8 that 
$$\mbox{G-dim\,gr}M\geq \mbox{G-dim}M\geq{\rm grade}M={\rm grade\,gr}M.$$
If ${\rm gr}M$ is perfect, then above inequalities are equalities. Hence $M\in{\mathcal G}_e$.

\medskip

(ii) Let $M\in{\mathcal G}_e$ with G-dim$M=d$. Then every syzygy $\Omega^iM$ of $M$ is also in ${\mathcal G}_e$. For, as ${\rm gr}(\Omega^iM)$ and $\Omega^i({\rm gr}M)$ are stably isomorphic, we see that G-dim\,$\Omega^iM=\mbox{G-dim\,gr}(\Omega^iM)={\rm max}\{0,\,d-i\}.$

\medskip

Applying Theorem 2.8 to the case that ${\rm gr}\Lambda$ is a Gorenstein ring, we get the following.

\medskip

2.10. {\sc Corollary.}\ \ 
{\sl Let $\Lambda$ be a filtered ring such that ${\rm gr}\Lambda$ is a commutative Gorenstein ring and $M$ a filtered $\Lambda$-module with a good filtration. Then the equality ${\rm grade}_{\Lambda}M={\rm grade}_{{\rm gr}\Lambda}{\rm gr}M$ holds. }

{\it Proof.}
Since all the finitely generated ${\rm gr}\Lambda$-modules have finite G-dimension (see the proof of \cite{AB}, Theorem 4.20), this follows from Theorem 2.8.   $\square$

\medskip

2.11. {\sc Theorem.}\ \ 
{\sl Let $\Lambda$ be a Gorenstein filtered ring. Let $M$ be a filtered $\Lambda$-module with a good filtration. Then the following equality holds.
$${\rm grade}M+\,^*{\rm dim\,gr}M=\,^*{\rm dim\,gr}\Lambda=\,^*{\rm id\,gr}\Lambda.$$  }

{\it Proof.}
This follows from A.9, A.10, A.12 and 2.8. $\square$

\medskip

When $\Lambda$ is a Gorenstein filtered ring, due to the above equality, we can define a holonomic module. Put $^*{\rm id\,gr}\Lambda=n$ and ${\rm id}\Lambda=d$. Let $M$ be a filtered $\Lambda$-module with a good filtration. Since ${\rm grade}M\leq{\rm id}\Lambda=d$, we have $n-\ ^*{\rm dim\,gr}M\leq d$, hence
$$^*{\rm dim\,gr}M\geq n-d.$$
This inequality is a generalization of Bernstein's inequality for a Weyl algebra (\cite{Bj1}). 

According to the case of Weyl algebras, we call a finitely generated filtered $\Lambda$-module $M$ a {\it holonomic module}, if $^*{\rm dim\,gr}M=n-d$. 

\vspace{1cm}

\begin{center}
3. {\sc Cohen-Macaulay modules and holonomic modules}
\end{center}

Throughout this section, we assume that $\Lambda$ is a filtered ring such that ${\rm gr}\Lambda$ is a Cohen-Macaulay $^*$local ring with the condition (P) (cf. Appendix). Let $M$ be a finitely generated filtered $\Lambda$-module such that $M\in {\mathcal G}$, i.e., G-dim\,gr$M<\infty$. It follows from 2.3, 2.8, A.8 and A.12 that the following holds:
$$(1)\ \ \ \mbox{\rm G-dim}M+\,^*{\rm depth\,gr}M\leq n$$
$$(2)\ \ \ {\rm grade}M+\,^*{\rm dim\,gr}M=n,$$
where we put $n:=\,^*{\rm depth\,gr}\Lambda=\,^*{\rm dim\,gr}\Lambda$. We say that $M\in\mathcal G$ is a {\it Cohen-Macaulay $\Lambda$-module of codimension $k$}, if $\,^*{\rm depth\,gr}M=\,^*{\rm dim\,gr}M=n-k$. Then it is easily seen that if $M$ is Cohen-Macaulay of codimension $k$ then it is perfect of grade $k$, where, due to \cite{AB}, Definition 4.34, we call $M$ perfect if G-dim$M={\rm grade}M$. Note also that $M$ is Cohen-Macaulay if and only if ${\rm gr}M$ is a perfect ${\rm gr}\Lambda$-module by A.8 and A.12. We put
$${\mathcal C}_k(\Lambda):=\{M\in{\mathcal G}:M\mbox{ is a Cohen-Macaulay $\Lambda$-module of codimension }k\}.$$
The following is an easy consequence of (1) and (2).

\medskip

3.1. {\sc Proposition.}\ \ 
{\sl Let $M\in{\mathcal C}_k(\Lambda)$. Then ${\rm Ext}_{\Lambda}^i(M,\Lambda)=0$ for all $i\neq k\ (i\geq 0).$}

\medskip

We slightly generalize \cite{I}, Lemma 2.7 and Theorem 2.8, and \cite{Iw}, as follows.

\medskip

3.2. {\sc Theorem.}\ \ 
{\sl Let $M\in{\mathcal G}$. 

\medskip

i) \ If $M\in{\mathcal C}_k(\Lambda)$, then ${\rm Ext}_{\Lambda}^k(M,\Lambda)\in{\mathcal C}_k(\Lambda^{\rm op}).$

\medskip

ii) The functor ${\rm Ext}_{\Lambda}^k(-,\Lambda)$ induces a duality between the categories $\underline{{\mathcal C}}_k(\Lambda)$ and $\underline{{\mathcal C}}_k(\Lambda^{\rm op}).$}

\medskip

3.2.1. {\sc Lemma.}\ \ 
{\sl Let $N$ be a finitely generated filtered $\Lambda$-module of ${\rm grade\,gr}N=s$. If the G-dimension of ${\rm gr}N$ is finite, then we have an embedding ${\rm gr}({\rm Ext}_{\Lambda}^s(N,\Lambda))\hookrightarrow{\rm Ext}_{{\rm gr}\Lambda}^s({\rm gr}N,{\rm gr}\Lambda).$ Moreover, if ${\rm gr}N$ is perfect, then the embedding is an isomorphism.}

{\it Proof.} 
Let $\cdots\rightarrow F_1\rightarrow F_0\rightarrow N\rightarrow0$ be a filtered free resolution of $N$. We use the notation of 2.8.1. It follows from 2.8.2 and 2.8 that 
$$E_s^1\cong H^s(E_{\bullet}^0)\cong H^{s-1}(F_{\bullet})={\rm Ext}_{{\rm gr}\Lambda}^{s-1}({\rm gr}N,{\rm gr}\Lambda)=0,$$
where a complex $F_{\bullet}: 0\rightarrow F_0^*\rightarrow F_1^* \rightarrow\cdots$ is as in 2.8.1. There exists an exact sequence 
$$0\rightarrow E_{s+1}^{r+1} \rightarrow E_{s+1}^r\rightarrow E_{s+2}^r$$
\noindent for all $r\geq1$ by 2.8.3, so that $E_{s+1}^r\subset E_{s+1}^1$ for all $r\geq1$. It follows from Lemma 2.8.5 that, for $r>\!>0$, 
$$E_{s+1}^r\cong{\rm gr}({\rm Ext}_{\Lambda}^s(N,\Lambda)).$$
Thus, by 2.8.2, we get
$${\rm gr}({\rm Ext}_{\Lambda}^s(N,\Lambda))\subset E_{s+1}^1\cong {\rm Ext}_{{\rm gr}\Lambda}^s({\rm gr}N,{\rm gr}\Lambda).$$
\noindent Assume further that ${\rm gr}N$ is perfect. Since $E_{s+2}^r$ is a subfactor of $E_{s+2}^1\cong {\rm Ext}_{{\rm gr}\Lambda}^{s+1}({\rm gr}N,{\rm gr}\Lambda)=0$, we see $E_{s+2}^r=0$, which shows that the embedding is an isomorphism.  $\square$

\medskip

3.2.2. {\sc Proof of 3.2.}\ \ 
i) Since gr$M$ is perfect of grade $k$, it holds that ${\rm Ext}_{{\rm gr}\Lambda}^k({\rm gr}M,{\rm gr}\Lambda)$ is perfect of grade $k$ by \cite{AB}, Proposition 4.35 and its proof, and so ${\rm gr\,Ext}_{\Lambda}^k(M,\Lambda)$ is perfect by Lemma 3.2.1. Hence ${\rm Ext}_{\Lambda}^k(M,\Lambda)\in{\mathcal C}_k(\Lambda^{\rm op})$. ii) Consider the exact sequence 
$$ 0\rightarrow {\rm Ext}_{\Lambda}^k(M,\Lambda)\rightarrow{\rm Tr}\Omega^{k-1}M\rightarrow\Omega{\rm Tr}\Omega^kM\rightarrow0$$
(see, for example, the proof of \cite{HN}, Lemma 2.1) and apply $(-)^*$ to it. Then we get a long exact sequence
$${\rm Ext}_{\Lambda^{\rm op}}^{k+1}({\rm Tr}\Omega^kM,\Lambda)\rightarrow{\rm Ext}_{\Lambda^{\rm op}}^k({\rm Tr}\Omega^{k-1}M,\Lambda)\rightarrow{\rm Ext}_{\Lambda^{\rm op}}^k({\rm Ext}_{\Lambda}^k(M,\Lambda),\Lambda)
\rightarrow {\rm Ext}_{\Lambda^{\rm op}}^{k+2}({\rm Tr}\Omega^kM,\Lambda).$$

\noindent Since G-dim${\rm Tr}\Omega^kM=0$ by assumption, the first and fourth terms of the above exact sequence vanishes. Hence $M\cong {\rm Ext}_{\Lambda^{\rm op}}^k({\rm Tr}\Omega^{k-1}M,\Lambda)\cong {\rm Ext}_{\Lambda^{\rm op}}^k({\rm Ext}_{\Lambda}^k(M,\Lambda),\Lambda)$ by \cite{HN}, Lemma 2.5. Therefore, there is a natural isomorphism $M\cong {\rm Ext}_{\Lambda^{\rm op}}^k({\rm Ext}_{\Lambda}^k(M,\Lambda),\Lambda)$ for $M\in{\mathcal C}_k(\Lambda)$, which induces a duality between the categories $\underline{{\mathcal C}}_k(\Lambda)$ and $\underline{{\mathcal C}}_k(\Lambda^{\rm op}). \ \ \square$

\medskip

3.2.3. {\sc Remark.}\ \ 
The proof 3.2.2 ii) only needs $M$ to be perfect with ${\rm grade}M=k$. Hence we see that if $M$ is perfect of grade $k$ then $M\cong {\rm Ext}_{\Lambda^{\rm op}}^k({\rm Ext}_{\Lambda}^k(M,\Lambda),\Lambda).$

\bigskip

We shall study holonomic modules when $\Lambda$ is a Gorenstein filtered ring, that is, ${\rm gr}\Lambda$ is Gorenstein, and generalize the former theory which is under the assumption of regularity (cf. \cite{HVO}, Chapter I\!I\!I, \S4). The assumption that $\Lambda$ is Gorenstein implies that ${\mathcal G}={\rm fil}\Lambda$ by A.9, where ${\rm fil}\Lambda$ is the category of all finitely generated filtered (left) $\Lambda$-modules. We recall from Corollary 2.4 and the end of section two that ${\rm id}_{\Lambda}\Lambda={\rm id}_{\Lambda^{\rm op}}\Lambda(=d)$ and $M\in{\rm fil}\Lambda$ is called holonomic, if $^*{\rm dim\,gr}M=n-d$, where $n=\,^{*}{\rm depth\,gr}\Lambda=\,^{*}{\rm dim\,gr}\Lambda=\,^{*}{\rm id\,gr}\Lambda$. We see that if $M\in {\mathcal C}_d(\Lambda)$ then $M$ is holonomic. We also note that $M$ is holonomic if and only if grade$M=d$ (or grade\,gr$M=d$) if and only if $M$ is perfect of grade $d$. We keep to assume $\Lambda$ to be a Gorenstein filtered ring and $d={\rm id}_{\Lambda}\Lambda$ in the rest of this section. According to \cite{Bj3}, Theorem 3.9, if $\Lambda$ is a Gorenstein filtered ring, then $\Lambda$ satisfies `Auslander condition' : 
\begin{center}
For every finitely generated $\Lambda$-module $M$ and integer $k\geq 0$, \\
it holds that 
${\rm grade}_{\Lambda^{\rm op}}N\geq k$ for all $\Lambda^{\rm op}$-submodules 
$N\subset {\rm Ext}_{\Lambda}^k(M,\Lambda)$.
\end{center}

\medskip

3.3. {\sc Proposition.}\ \ 
{\sl Let $M$ be a finitely generated filtered $\Lambda$-module. Let $M$ be holonomic and $N$ a $\Lambda$-submodule of $M$. Then $N,\ M/N$ are holonomic.}

{\it Proof.}
It follows from \cite{HN}, Lemma 2.11 (cf. also \cite{Bj3}, Theorem 3.9) that ${\rm grade}N\geq d$ and ${\rm grade}M/N\geq d$, so that ${\rm grade}N=d$ and ${\rm grade}M/N=d.\ \ \square$

\medskip

3.4. {\sc Proposition.}\ \ 
{\sl A holonomic module is artinian. Therefore, it is of finite length.}

\medskip

We use the following easy lemma for a proof.

\medskip

3.4.1. {\sc Lemma. \ \ }
{\sl Let $M_i\ (i=0,1,\cdots)$ be a module over a ring and $f_i:M_i\rightarrow M_{i+1}\ (i=0,1,\cdots)$ is a homomorphism. Assume that $M_0$ is Noetherian and $f_i\ (i=0,1,\cdots)$ is surjective. Then there exists an interger $m$ such that $f_i$ is an isomorphism for all $i\geq m$.}

\medskip

3.4.2. {\sc Proof of 3.4.}\ \ 
Let $M$ be a holonomic $\Lambda$-module and $M=M_0\supset M_1\supset\cdots$ a descending chain of $\Lambda$-submodules of $M$. Then $M_i,\ M_{i-1}/M_i$ are holonomic ($i\geq 1$), and so, from an exact sequence $0\rightarrow M_i\rightarrow M_{i-1}\rightarrow M_{i-1}/M_i\rightarrow0$, we get an exact sequence
$$0\rightarrow {\Bbb E}(M_{i-1}/M_i)\rightarrow {\Bbb E}M_{i-1}\rightarrow{\Bbb E}M_i\rightarrow 0,$$
where we put ${\Bbb E}(-)={\rm Ext}_{\Lambda}^d(-,\Lambda)$. By Lemma 3.4.1, there exists an integer $m$ such that ${\Bbb E}M_{i-1}\rightarrow{\Bbb E}M_i$ is an isomorphism for $i\geq m+1$. Hence ${\Bbb E}(M_{i-1}/M_i)=0$ for $i\geq m+1$. Hence $M_{i-1}/M_i=0$ for $i\geq m+1$ by Remark 3.2.3, that is, $M_m=M_{m+1}=\cdots$. This completes the proof.  $\square$

\medskip

We generalize \cite{HVO}, Chapter I\!I\!I, 4.2.18 Theorem (p. 194), which characterizes a holonomic module by its associated graded module. We put ${\rm Min}({\rm gr}M)=\{{\frak p}:{\frak p}$ is a minimal element of ${\rm Supp}({\rm gr}M)\}$ for $M\in {\rm fil}\Lambda$.

\medskip

3.5. {\sc Theorem.}\ \ 
{\sl Let $M\in{\rm fil}\Lambda.$ Then the following are equivalent.

\medskip

(1) $M$ is holonomic,

\medskip

(2) ${\rm ht}{\frak p}=d$ for all ${\frak p}\in {\rm Min}({\rm gr}M).$}

\medskip

A finitely generated module $M$ over a two-sided Noetherian ring is called {\it pure}, if ${\rm grade}N={\rm grade}M$ for all nonzero submodules $N$ of $M$. 

\medskip

3.5.1. {\sc Lemma.}\ \ 
{\sl Let $M\in {\rm fil}\Lambda.$ Then $M$ is pure if and only if ${\rm gr}M$ is a pure ${\rm gr}\Lambda$-module under a suitable filtration on $M$.}

{\it Proof.}
Let $M$ be pure. Put $s={\rm grade}M$ and $N:={\rm Ext}_{\Lambda}^s(M,\Lambda)$. Since $\Lambda$ satisfies Auslander condition, it follows that ${\rm grade}N=s$ by \cite{HN}, Lemma 2.8, so that ${\rm grade\,gr}N=s$ by 2.8, hence ${\rm Ext}_{{\rm gr}\Lambda}^s({\rm gr}N,{\rm gr}\Lambda)$ is pure by \cite{HN}, Proposition 2.13. By 3.2.1, we have ${\rm gr\,Ext}_{\Lambda^{\rm op}}^s(N,\Lambda)\subset{\rm Ext}_{{\rm gr}\Lambda}^s({\rm gr}N,{\rm gr}\Lambda).$ Hence ${\rm gr\,Ext}_{\Lambda^{\rm op}}^s(N,\Lambda)$ is a pure ${\rm gr}\Lambda$-module. By \cite{HN}, Theorem 2.3, there exists an exact sequence
$$0\rightarrow{\rm Ext}_{\Lambda^{\rm op}}^{s+1}({\rm Tr}\Omega^{s}M,\Lambda)\rightarrow M\rightarrow{\rm Ext}_{\Lambda^{\rm op}}^s(N,\Lambda).$$
Since ${\rm grade\,Ext}_{\Lambda^{\rm op}}^{s+1}({\rm Tr}\Omega^{s}M,\Lambda)\geq s+1$ by Auslander condition and $M$ is pure, we see ${\rm Ext}_{\Lambda^{\rm op}}^{s+1}({\rm Tr}\Omega^{s}M,\Lambda)=0$. Therefore, $M\subset {\rm Ext}_{\Lambda^{\rm op}}^s(N,\Lambda)$. According to a filtration on $M$ induced from that of ${\rm Ext}_{\Lambda^{\rm op}}^s(N,\Lambda)$, we get an inclusion ${\rm gr}M\subset{\rm gr\,Ext}_{\Lambda^{\rm op}}^s(N,\Lambda)$, hence ${\rm gr}M$ is pure. The converse is obvious by Theorem 2.8.   $\square$

\medskip

3.5.2. {\sc Lemma.}\ \ 
{\sl Let $R$ be a commutative Gorenstein ring and $M'$ a pure $R$-module. Then ${\rm grade}M'={\rm dim} R_{\frak p}$ for each ${\frak p}\in {\rm Min}(M').$}

{\it Proof.}
Since $R_{\frak p}$ is a Gorenstein local ring, we have an equality ${\rm grade}M^{'}_{\frak p}+{\rm dim}M^{'}_{\frak p}={\rm dim}R_{\frak p}$ (cf. \cite{GN}, Proposition 4.11). Since ${\frak p}$ is minimal, ${\rm dim}M^{'}_{\frak p}=0$, so that, ${\rm grade}M^{'}_{\frak p}={\rm dim}R_{\frak p}$. 

Put $g={\rm grade}M^{'}_{\frak p},\ g'={\rm grade}M'$. Since ${\rm Ext}_R^g(M',R)_{\frak p}={\rm Ext}_{R_{\frak p}}^g(M_{\frak p}^{'},R_{\frak p})\neq0$, we have ${\rm Ext}_R^g(M',R)\neq0$. Hence $g\geq g'$ holds. Suppose that ${\rm Ext}_R^k({\rm Ext}_R^k(M',R),R)\neq0$ for $k>g'$. Then there exists $N\subset M'$ such that ${\rm grade}N=k>g'$ by \cite{HN}, Theorem 2.3, which contradicts the purity of $M'$. Hence ${\rm Ext}_R^k({\rm Ext}_R^k(M',R),R)=0$ for all $k>g'$. But by A.15, ${\rm grade\,Ext}_{R_{\frak p}}^g(M_{\frak p}^{'},R_{\frak p})=g$. Therefore, we see $g\leq g'$, and so, $g=g'$. This completes the proof.   $\square$

\medskip

3.5.3. {\sc proof of Theorem 3.5.}\ \ 
Put $R={\rm gr}\Lambda$.

(1)$\Rightarrow$(2): Assume that $M$ is holonomic. Since $M$ is pure by Proposition 3.3, ${\rm gr}M$ is pure by 3.5.1. Thus $d={\rm grade\,gr}M={\rm dim}R_{\frak p}$ for all ${\frak p}\in{\rm Min}({\rm gr}M)$ by 3.5.2. Therefore, ${\rm ht}{\frak p}=d$ for all ${\frak p}\in{\rm Min}({\rm gr}M)$.

(2)$\Rightarrow$(1): Put $I=[0:_R{\rm gr}M]$. Since $R$ is Cohen-Macaulay, we have ${\rm ht}I={\rm grade}R/I$ by \cite{BH}, Corollary 2.1.4. It follows from \cite{BH}, Proposition 1.2.10(e) that ${\rm grade}R/I={\rm grade\,gr}M$. By assumption, ${\rm ht}I=d$, so that, ${\rm grade\,gr}M=d$, that is, ${\rm grade}M=d$ by 2.8. Hence $M$ is holonomic.  $\square$

\medskip

A module having higher grade has a good property.

\medskip

3.6. {\sc Proposition}.\ \ 
{\sl Let $M$ be a finitely generated filtered $\Lambda$-module with ${\rm grade}M=\ell$, where $\ell=d-1\mbox{ \sl or }d-2$. Then $M$ is a perfect $\Lambda$-module if and only if there exists a finitely generated filtered $\Lambda^{\rm op}$-module $M'$ of grade $\ell$ with $M\cong {\rm Ext}_{\Lambda^{\rm op}}^{\ell}(M',\Lambda)$. }

{\it Proof.}
Assume that $M\cong {\rm Ext}_{\Lambda^{\rm op}}^{\ell}(M',\Lambda)$ with ${\rm grade}M'=\ell$. 

The case $\ell=d-1$: We see ${\rm grade}M=d-1$ by assumption. It follows that ${\rm grade}\,{\rm Ext}_{\Lambda}^d(M,\Lambda)={\rm grade}\,{\rm Ext}_{\Lambda}^d({\rm Ext}_{\Lambda}^{d-1}(M',\Lambda),\Lambda)\geq d+2$ by \cite{HN}, Corollary 2.10. This shows that ${\rm Ext}_{\Lambda}^d(M,\Lambda)=0$, that is, ${\rm G\mbox{-}dim}M\leq d-1$. Hence ${\rm G\mbox{-}dim}M={\rm grade}M=d-1$, so that $M$ is perfect.

The case $\ell=d-2$: It follows from the similar computations as the above case that ${\rm grade}\,{\rm Ext}_{\Lambda}^d(M,\Lambda)\geq d+2$ and ${\rm grade}\,{\rm Ext}_{\Lambda}^{d-1}(M,\Lambda)\geq d+1$. Hense ${\rm Ext}_{\Lambda}^d(M,\Lambda)={\rm Ext}_{\Lambda}^{d-1}(M,\Lambda)=0$, so that ${\rm G\mbox{-}dim} M={\rm grade}M=d-2$, i.e., $M$ is perfect.

The converse follows from \cite{Iw}, Theorem 4.  $\square$

\medskip

3.7. \ Following \cite{HVO}, Chapter I\!I\!I, 4.3, we call a filtered $\Lambda$-module $M$ {\it geometrically pure} ({\it geo-pure} for short), if ${\rm dim}_{{\rm gr}\Lambda}{\rm gr}M={\rm dim}({\rm gr}\Lambda/{\frak p})$ for all ${\frak p}\in {\rm Min}({\rm gr}M)$. Then we have the following proposition which is a generalization of \cite{HVO}, Chapter I\!I\!I, 4.3.6 Corollary.

\medskip

3.7.1. {\sc Proposition}.\ \ 
{\sl Let $M$ be a finitely generated filtered $\Lambda$-module, and put 
${\rm Assh\,gr}M:=\{\,\frak{p}\in {\rm Supp\,gr}M \ | \ {\rm dim}\,{\rm gr}\Lambda/\frak{p}={\rm dim\,gr}M\,\}$. 
Then the following conditions are equivalent.

\medskip

(1) $M$ is pure,

\medskip

(2) $M$ is geo-pure and ${\rm gr}M$ has no embedded prime.}

\medskip

(3) ${\rm Ass\,gr}M={\rm Assh\,gr}M$

{\it Proof.} 
(1)$\Rightarrow$(2): Let $M$ be pure. Then ${\rm gr}M$ is pure by 3.5.1. Take any ${\frak p}\in{\rm Min}({\rm gr}M)$. Since ${\frak p}\in{\rm Ass\,gr}M$, we have ${\rm gr}\Lambda/{\frak p}\hookrightarrow{\rm gr}M$, so ${\rm grade\,gr}\Lambda/{\frak p}={\rm grade\,gr}M$. Using Theorem A.12, we have ${\rm dim\,gr}\Lambda/{\frak p}={\rm dim\,gr}M$. Hence $M$ is geo-pure. Take any ${\frak p}\in{\rm Ass\,gr}M$, then ${\rm gr}\Lambda/{\frak p}\hookrightarrow{\rm gr}M$. Thus ${\rm dim\,gr}\Lambda/{\frak p}={\rm dim\,gr}M$, by A.12. Therefore, ${\rm Ass\,gr}M={\rm Min\,gr}M$, i.e., ${\rm gr}M$ has no embedded primes.

(2)$\Rightarrow$(3): The former condition implies ${\rm Assh\,gr}M={\rm Min\,gr}M$, 
and the latter one implies ${\rm Min\,gr}M={\rm Ass\,gr}M.$

(3)$\Rightarrow$(1): By 3.5.1, it suffices to prove that ${\rm gr}M$ is pure. Let $N$ be a 
${\rm gr}\Lambda$-submodule of ${\rm gr}M$. Take any ${\frak p}\in{\rm Ass}N$. Then 
${\rm gr}\Lambda/{\frak p}\hookrightarrow N$. Thus, by A.12 and assumption, we have 
${\rm grade\,gr}\Lambda/{\frak p}={\rm grade\,gr}M$. By \cite{HN}, Lemma 2.11, we have 
$${\rm grade\,gr}M\leq{\rm grade}N\leq{\rm grade\,gr}\Lambda/{\frak p}={\rm grade\,gr}M.$$
Hence ${\rm grade\,gr}M={\rm grade}N$. This completes the proof.    $\square$  

\medskip

3.8. {\sc Example}.\ \ 
We provide an example of a Gorenstein filtered ring $\Lambda$.

\noindent Let $R=k[[x^2,x^3]]$ be a subring of a formal power series ring $k[[x]]$, where $k$ is a field of characteristic zero. Then $(R,{\frak m})$ is a local Gorenstein (non-regular) ring of dim$R=1$, where ${\frak m}=(x^2,x^3)$. Let a differential operator $T=x\partial$ with $\partial=d/dx$. Let $\Lambda$ be a subring of the first Weyl algebra (see \cite{Bj1}, \cite{HVO}) generated by $R$ and $T$. Then every element of $\Lambda$ is written as $\Sigma a_iT^i,\ a_i\in R$. Note that $Tx^i=x^iT+ix^i,\ i\geq 2$. For $P=\Sigma a_iT^i\in\Lambda$, we put ${\rm ord}P={\rm max}\{i:a_i\neq 0\}$, an order of $P$. Let ${\mathcal F}_i\Lambda:=\{P\in\Lambda:{\rm ord}P\leq i\}$. Then $\{{\mathcal F}_i\Lambda\}$ is a filtration of $\Lambda$ and ${\rm gr}\Lambda=R[t]$, where $t=\sigma_1(T)$. Thus ${\rm gr}\Lambda$ is Gorenstein $^*$local of dimension 2. Note that ${\frak m}+tR[t]$ is a unique $^*$maximal ideal.

\medskip

1) ${\rm id}\Lambda=2$

\noindent Let $I:=\Lambda T+\Lambda x^2$ be a left ideal of $\Lambda$. Then $I\neq\Lambda$. We put induced filtrations to $I$ and $\Lambda/I$. , i.e., 
$${\mathcal F}_iI=I\cap{\mathcal F}_i\Lambda,\ \ {\mathcal F}_i(\Lambda/I)=({\mathcal F}_i\Lambda+I)/I,\ \ i\geq 0.$$
Then $0\rightarrow I\rightarrow\Lambda\rightarrow\Lambda/I\rightarrow 0$ is a strict exact sequence. Hence $0\rightarrow {\rm gr}I\rightarrow{\rm gr}\Lambda\rightarrow{\rm gr}(\Lambda/I)\rightarrow 0$ is exact. Since gr$I$ contains $t$ and $x^2$, ${\rm gr}(\Lambda/I)={\rm gr}\Lambda/{\rm gr}I$ is an Artinian gr$\Lambda$-module. Hence dim$_{{\rm gr}\Lambda}{\rm gr}(\Lambda/I)=0$. Thus grade\,$\Lambda/I=2$ by Theorem 2.11, and then id$\Lambda=2$ by Corollary 2.4. 
So $\Lambda/I$ is holonomic. 

\medskip

2) ${\rm gl\,dim}\Lambda=\infty$

\noindent It follows that ${\rm pd}_R(R/{\frak m})$ is infinite. Since $\Lambda$ is $R$-free, we see that ${\rm pd}_{\Lambda}(\Lambda/\Lambda{\frak m})$ is also infinite (cf. \cite{R}, Exercise 9.21, p.258)

\medskip 

\vspace{1cm}

\begin{center}
{\sc Appendix}
\end{center}

In Appendix, we provide the fact about graded rings, especially $^*$local rings .

\medskip

\begin{center}
1. {\sc Summary for $^*$local rings}\ \  
\end{center}

Let $R$ be a commutative Noetherian ring. We gather some facts about a graded ring. For the detail, the reader is referred to \cite{BH}, \cite{GW}, and \cite{NVO}.

A ring $R$ is called a {\it graded ring}, if 

\medskip

i)  $R=\oplus_{i\in \Bbb Z}R_i$ as an additive group,

\medskip

ii) $R_iR_j\subset R_{i+j}$ for all $i,j\in {\Bbb Z}$.

\medskip

An $R$-module $M$ is called a {\it graded module}, if 

\medskip

i)  $M=\oplus_{i\in\Bbb Z}M_i$ as an additive groups,

\medskip

ii) $R_iM_j\subset M_{i+j}$ for all $i,\ j\in{\Bbb Z}$.

\medskip

An $R$-homomorphism $f:M \rightarrow N$ of graded modules is called a {\it graded homomorphism}, if $f(M_i)\subset N_i$ for all $i\in {\Bbb Z}$. All graded modules in  ${\rm mod}R$ and all graded homomorphisms form the category of graded modules, which we denote by ${\rm mod}_0R$.

A graded submodule of a graded ring $R$ is called a {\it graded ideal}. For any ideal $I$ of $R$, we denote by $I^*$ the graded ideal generated by all homogeneous elements of $I$. A graded ideal ${\frak m}$ of $R$ is called $^*${\it maximal}, if it is a maximal element of all proper graded ideals of $R$. We say that $R$ is a $^*${\it local} ring, if $R$ has a unique $^*$maximal ideal ${\frak m}$. A $^*$local ring $R$ with the $^*$maximal ideal ${\frak m}$ is denoted by $(R,{\frak m})$. The theory of $^*$local ring is well developed and a lot of facts that hold for local rings also hold for $^*$local rings (see \cite{BH} and \cite{GW}).  

Let $M$ be a finite $R$-module. For an ideal $I$, we denote $I$-depth of $M$ by  ${\rm depth}(I,M)$(\cite{M}). Let $(R,{\frak m})$ be a $^*$local ring and $M\in {\rm mod}R$. We put $^*{\rm depth}M:={\rm depth}({\frak m},M)$. We shall use $^*$depth as a substitute of depth for a local ring. 

A graded module $M$ over a graded ring $R$ is called a $^*${\it injective} module, if it is an injective object in ${\rm mod}_0R$(\cite{BH}, \S 3.6). We denote by $^*{\rm id}M$ the $^*$injective dimension of $M$. By definition, $^*{\rm id}M\leq k$ if and only if there exists a minimal $^*$injective resolution 
$$0\rightarrow M \rightarrow \:^*\negmedspace E^0(M)\rightarrow\cdots \rightarrow \:^* \negmedspace E^k(M)\rightarrow 0.$$  
It is easily seen that $^*{\rm id}M\leq k$ if and only if ${\rm Ext}_R^i(N,M)=0$ for all $i>k$ and all $N\in{\rm mod}_0R$.

Let $(R,{\frak m})$ be a $^*$local ring. Consider the following condition.
\begin{center}
(P) There exists an element of positive degree in $R-{\frak p}$ \\ 
for any graded prime ideal ${\frak p\neq \frak m}$
\end{center}

A positively graded ring satisfies the condition (P). The other examples are seen in \cite{NVO}, Chapter B, I\!I\!I, 3.2.

The following is known.

\medskip

A.1. {\sc Proposition}.\ \
{\sl Let $(R,\frak m)$ be a $^*$local ring with the condition (P). Then, for every graded ideal $\frak a$ and every set of graded prime ideals $\frak p_1,\cdots,\frak p_n$, there exists $i$ such that $\frak a\subset \frak p_i$, whenever all homogeneous elements of $\frak a$ are contained in $\cup_{i=1}^n\frak p_i$.}

{\it Proof.} See \cite{MR}, Lemma 2. $\square$

\medskip

Using Proposition A.1, the following is proved as the local case.

\medskip

A.2. {\sc Proposition}.\ \ 
{\sl Let $(R,\frak m)$ be a $^*$local ring with the condition (P). Let $M$ be a finite graded module with $^*{\rm depth}M=t$. Then there exists an $M$-sequence $x_1,\cdots,x_t$ consisting of homogeneous elements in $\frak m$.}

\medskip

We note the following graded version of Nakayama's Lemma.

\medskip

A.3. {\sc Lemma}.\ \ 
{\sl Let $(R,\frak m)$ be a $^*$local ring and $M$ a finite graded $R$-module. If $\frak m M=M$, then $M=0$.}

\medskip

In the following, we assume that $(R,\frak m)$ is a $^*$local ring with the condition (P).

\medskip
A.4. {\sc Lemma}.\ \ 
{\sl Let $M,\ N$ be the non-zero finite graded $R$-module with $^*{\rm depth}N=0$. Then ${\rm Hom}_R(M,N)\neq 0$.}

{\it Proof.} It is well-known, so we omit the proof. $\square$

\medskip

A.5. {\sc Corollary}.\ \ 
{\sl Assume that $^*{\rm depth}R=0$. Let $M$ be a finite graded $R$-module. Then $M^*=0$ implies $M=0$.}

\medskip

We state the graded version of \cite{AB}, 4.11-13 in the following A.6-A.8.

\medskip

A.6. {\sc Proposition}.\ \ 
{\sl Assume that $^*{\rm depth}R=0$. Let $M$ be a finite graded $R$-module. Then {\rm G}-${\rm dim}M<\infty$ if and only if {\rm G}-${\rm dim}M=0$.}

{\it Proof.} It suffices to prove that G-dim$M<\infty$ implies G-dim$M=0$. 

Suppose that G-dim$M\leq 1$. We have an exact sequence $0\rightarrow L_1\rightarrow L_0\rightarrow M\rightarrow 0$ with G-dim$L_i=0\ (i=0,1)$. Hence we have an exact sequence
$$0\rightarrow M^*\rightarrow L_0^*\rightarrow L_1^*\rightarrow{\rm Ext}_R^1(M,R)\rightarrow 0$$
and ${\rm Ext}_R^i(M,R)=0$ for $i>1$. By this sequence, we have an exact sequence
$$0\rightarrow{\rm Ext}_R^1(M,R)^*\rightarrow L_1\rightarrow L_0,$$
where $L_1\rightarrow L_0$ is monic. Thus ${\rm Ext}_R^1(M,R)^*=0$, and so ${\rm Ext}_R^1(M,R)=0$ by A.5 Corollary.

Suppose that G-dim$M\leq n$.  Let $0\rightarrow L_n\overset{f_n}{\rightarrow}\cdots \overset{f_1}{\rightarrow} L_0\rightarrow M\rightarrow 0$ be exact with G-dim$L_i=0\ (0\leq i\leq n)$. Since G-dim(Im$f_{n-1})\leq 1$, we have G-dim(Im$f_{n-1})=0$ by the above argument. Repeating this process, we get G-dim$M=0$. $\square$

\medskip

We want to generalize \cite{AB}, Theorem 4.13 (b) to the graded case. The proof of it needs a part of \cite{AB}, Proposition 4.12. Thus we adapt this proposition as follows.
  
\medskip

A.7. {\sc Proposition}.\ \ 
{\sl Assume that $^*{\rm depth}R=t$. Let $M$ be a finite graded $R$-module with {\rm G}-${\rm dim}M<\infty$. Then the following are equivalent.

\medskip
(1) {\rm G}-${\rm dim}M=0$.

\medskip

(2) $^*{\rm depth}M\geq \  ^*{\rm depth}R$.

\medskip

(3) $^*{\rm depth}M= \ ^*{\rm depth}R$.}

\medskip

{\it Proof.} (1) $\Rightarrow$ (2): Let $x_1,\cdots,x_i$ be a homogeneous regular sequence in ${\frak m}$. We show that $x_1,\cdots,x_i$ is an $M$-sequence by induction on $i$. Let $i=1$. Since $M\cong M^{**}$ is torsionfree, $x_1$ is $M$-regular.

Suppose that $i>1$ and the assertion holds for $i-1$. Then $x_1,\cdots,x_{i-1}$ is an $M$-sequence. Put $I=(x_1,\cdots,x_{i-1}),\ \overline{R}=R/I,\ \overline{M}=M/IM$. Then $(\overline{R},{\frak m}/I)$ is a $^*$local ring with the condition (P). By \cite{AB}, Lemma 4.9, G-dim$_{\overline{R}}\overline{M}=\mbox{G-dim}_RM=0$. Since $\overline{x}_i\in\overline{R}$ is a regular element, $\overline{x}_i$ is $\overline{M}$-regular, hence $x_1,\cdots,x_i$ is an $M$-sequence. Therefore, $^*{\rm depth}M\geq\,^*{\rm depth}R$. 

(2) $\Rightarrow$ (1): By assumption, it suffices to prove that ${\rm Ext}_R^i(M,R)=0$ for $i>0$. We show the assertion by induction on $t=\,^*{\rm depth}R$.

Let $t=0$. Then G-dim$M=0$ by Proposition A.6

Let $t>0$. Then $^*{\rm depth}M\geq\,^*{\rm depth}R\geq1$. We take a homogeneous element $x\in{\frak m}$ which is $R$ and $M$-regular. Then, by \cite{BH}, 1.2.10 (d), 
$$^*{\rm depth}_{R/xR}M/xM=\,^*{\rm depth}_RM-1\geq\,^*{\rm depth}R-1=\,^*{\rm depth}R/xR.$$
Hence we have ${\rm Ext}_{R/xR}^i(M/xM,R/xR)=0$ for $i>0$ by induction. This gives ${\rm Ext}_R^i(M,R/xR)=0$ for $i>0$. From an exact sequence $0\rightarrow R\overset{x}{\rightarrow} R\rightarrow R/xR\rightarrow 0$, we get an exact sequence
$${\rm Ext}_R^i(M,R)\overset{x}{\rightarrow} {\rm Ext}_R^i(M,R)\rightarrow{\rm Ext}_R^i(M,R/xR)=0.$$
By Nakayama's Lemma, it holds that ${\rm Ext}_R^i(M,R)=0$ for $i>0$.  

(2) $\Rightarrow$ (3): When $^*{\rm depth}R=0$, we have G-dim$M=0$ by Proposition A.6. Since $^*{\rm depth}R=0$, we have an exact sequence $0\rightarrow R/{\frak m} \rightarrow R$ which gives an exact sequence
$$0\rightarrow {\rm Hom}_R(M^*,R/{\frak m})\rightarrow M^{**}\cong M.$$
Since $M^*\neq 0$, we have ${\rm Hom}_R(M^*,R/{\frak m})\neq 0$. Since ${\frak m}{\rm Hom}_R(M^*,R/{\frak m})=0$, we see that ${\frak m}$ has no $M$-regular element, so that $^*{\rm depth}M=0$. Thus (3) holds. 

Let $^*{\rm depth}R>0$. We have $^*{\rm depth}M\geq\,^*{\rm depth}R\geq 1$, so that there is a homogeneous element $x\in{\frak m}$ which is $R$ and $M$-regular. By \cite{AB}, Lemma 4.9, we have G-dim$_{R/xR}M/xM<\infty$. We have
$$^*{\rm depth}_{R/xR}M/xM=\,^*{\rm depth}_RM-1\geq\,^*{\rm depth}R-1=\,^*{\rm depth}R/xR.$$
Hence, by induction on $^*{\rm depth}R$, we have $^*{\rm depth}_{R/xR}M/xM=\,^*{\rm depth}R/xR,$
and then $^*{\rm depth}M=\,^*{\rm depth}R$. 

Since (3) $\Rightarrow$ (2) is obvious, we accomplish the proof. $\square$
 
\medskip

A.8. {\sc Theorem}.\ \ 
{\sl Let $M$ be a finite graded $R$-module with {\rm G}-${\rm dim}M<\infty$. Then we have an equality
$$\mbox{\rm G-dim}M+\ ^*{\rm depth}M=\ ^*{\rm depth}R$$ }
{\it Proof.} We state the proof which is an adaptation of \cite{AB}. If G-dim$M=0$, we are done by the previous proposition. Suppose that G-dim$M=n>0$ and the equation holds for $n-1$. Let $0\rightarrow K\rightarrow F\rightarrow M\rightarrow0$ be exact with $F$ graded free and $K$ a graded module. Since G-dim$K=n-1$, we have $\mbox{G-dim}K+\,^*{\rm depth}K=\,^*{\rm depth}R$ by induction. Suppose that $^*{\rm depth}M\geq\,^*{\rm depth}F=\,^*{\rm depth}R$. Then G-dim$M=0$ holds by the previous proposition. This contradicts to G-dim$M>0$. Hence $^*{\rm depth}M<\,^*{\rm depth}F$, so $^*{\rm depth}K=\,^*{\rm depth}M+1$ by, e.g., \cite{BH}, 1.2.9. Therefore, $n+\,^*{\rm depth}M=\,^*{\rm depth}R$. $\square$

\medskip

Let $M$ be a finite graded $R$-module. Then the similar argument to \cite{AB}, 4.14 and 4.15 shows that G-${\rm dim}M\leq n$ if and only if G-${\rm dim}M_{\frak p}\leq n$ for all graded prime (respectively, graded maximal) ideals ${\frak p}$ of $R$. Note that all the prime ideals in ${\rm Ass}M$ are graded ideals (e.g. \cite{BH}, Lemma 1.5.6). Thus, in $^*$local case, we have that G-${\rm dim}M\leq n$ if and only if G-${\rm dim}M_{\frak m}\leq n$. Thus we give the following characterization of Gorensteiness.

\medskip

A.9. {\sc Theorem}.\ \ 
{\sl  Let $(R,{\frak m})$ be a $^*$local ring with the condition (P). Then the following are equivalent.

\medskip

(1) $R$ is Gorenstein.

\medskip

(2) Every finite graded $R$-module has finite G-dimension.

\medskip

Under these equivalent conditions, the equality $^*{\rm id}R=\ ^*{\rm depth}R$ holds.}

\medskip

{\it Proof.} $(1)\Rightarrow (2)$: Since $R_{\frak m}$ is Gorenstein, we have G-${\rm dim}M_{\frak m}<\infty$, hence G-${\rm dim}M<\infty$ by above.

$(2)\Rightarrow(1)$: Let $t=\ ^*{\rm depth}R$. Take any finite graded $R$-module $M$. Since G-${\rm dim}M=t-\ ^*{\rm depth}M\leq t$ by Theorem A.8, we have that ${\rm Ext}_R^i(M,R)=0$ for all $i>t$. Hence $^*{\rm id}R\leq t$. It holds from \cite{BH}, Theorem 3.6.5 or \cite{NVO}, Chapter B, III.1.7 that ${\rm id}R\leq\ ^*{\rm id}R+1\leq t+1$. Hence $R$ is Gorenstein.

The second statement follows from the similar argument to the local case (cf. \cite{BH},Theorem 3.1.17). We note that `the residue field' in the local case should be replaced by `the unique graded simple module $R/{\frak m}$' in $^*$local case and the use of the graded version of Bass's Lemma (see e.g. \cite{NVO}, Chapter B, III.1.9) is effective. $\square$

\bigskip

Let $(R,{\frak m})$ be a $^*$local ring. Then one of the following cases occurs (\cite{GW}, \S1 or \cite{BH}, \S1.5):

\medskip

A. $R/{\frak m}$ is a field,

\medskip

B. $R/{\frak m}\cong k[t,t^{-1}]$, where $k$ is a field and $t$ is a homogeneous element of positive degree and transcendental over $k$.

\medskip

We put $^*$dim$R:={\rm ht}{\frak m}$ the $^*$dimension of a $^*$local ring $(R,{\frak m})$. Note that $^*$dim$R$ equals the supremum of all numbers $h$ such that there exists a chain of graded prime ideals ${\frak p}_0\subset{\frak p}_1\subset\cdots\subset{\frak p}_h$ in $R$ \cite{BH}. Let $M$ be a finite graded $R$-module. It is easily seen that $[0:_RM]$ is a graded ideal. Thus we put $^*$dim$M:=\ ^*{\rm dim}R/[0:_RM]$.

\medskip

A.10. {\sc Lemma}.\ \ 
{\sl Let $(R,{\frak m})$ be a Cohen-Macaulay $^*$local ring with the condition (P) and ${\rm dim}R=n$, and $M$ a finite graded $R$-module. Then we have }
$$^*{\rm dim}R=\ ^*{\rm depth}R=\left\{
\begin{array}{lcc} n &  \mbox{\sl for} & \mbox{\sl Case {\rm A},}\\
                    n-1 & \mbox{\sl for} & \mbox{\sl Case {\rm B}}.
\end{array}\right.$$ 
 
$$^*{\rm dim}M=\left\{
\begin{array}{lcc} {\rm dim}M &  \mbox{\sl for} & \mbox{\sl Case {\rm A},}\\
                    {\rm dim}M-1 & \mbox{\sl for} & \mbox{\sl Case {\rm B}}.
\end{array}\right.$$                   
\noindent {\sl Moreover, assume that $R$ is Gorenstein, then ${\rm id}R={\rm dim}R=n$, where {\rm id}$R$ stands for the injective dimension of $R$.}

{\it Proof.} Case A. Let $\frak n$ be a maximal ideal with ${\rm ht}\ {\frak n}=n$. If $\frak n=\frak m$, then ${\rm ht}\ {\frak m}=n$. Suppose that $\frak n$ is not equal to $\frak m$. Then $\frak n$ is not graded, so ${\rm ht}\ {\frak n}/{\frak n}^*=1$. Since $R_{\frak n}$ is Cohen-Macaulay, 
$${\rm ht}\ {\frak n}^*R_{\frak n}+{\rm dim}R_{\frak n}/{\frak n}^*R_{\frak n}={\rm dim} R_{\frak n}=n$$
 (\cite{M}, Theorem 17.4). Hence ${\rm ht}\ {\frak n}^*R_{\frak n}=n-1$, so ${\rm ht}\ {\frak n}^*=n-1$. Thus ${\rm ht}\ {\frak m}\geq {\rm ht}\ {\frak n}^*+1=n$, so that ${\rm ht}\ {\frak m}=n$. Therefore,
 $$^*{\rm depth}R={\rm depth}R_{\frak m}={\rm dim}R_{\frak m}={\rm ht}\ {\frak m}=n.$$
 
 Case B. Let $\frak n$ be the same as in Case A. Since $\frak n$ is not graded, we have ${\rm ht}\ {\frak n}^*=n-1$ by the similar way to Case A. By assumption, we have that ${\frak m}\supset{\frak n}^*$ and $\frak m$ is not maximal, so ${\frak m}={\frak n}^*$. Therefore, ${\rm ht}\ {\frak m}=n-1$, hence we get $^*{\rm depth}R=n-1$ by the similar way to Case A.
 
The equality concerning $^*$dim$M$ follows from the fact that cases A and B are preserved modulo $[0:_RM]$.
 
The latter statement is proved in \cite{B} more generally. $\square$
 
\medskip

A.11. {\sc Lemma}\ \ 
{\sl Let $(R,{\frak m})$ be a Cohen-Macaulay $^*$local ring with the condition (P) and $x$ a homogeneous element in $\frak m$. If $x$ is regular, then ${\rm dim}R/xR={\rm dim}R-1$.}

\medskip

{\it Proof.} The well-known induction argument works due to A.10 Lemma. $\square$

\medskip

A.12. {\sc Theorem}\ \ 
{\sl Let $(R,{\frak m})$ be a Cohen-Macaulay $^*$local ring with the condition (P) and $M$ a finite graded $R$-module. Then
$${\rm grade}M+{\rm dim}M={\rm dim}R$$ }
{\it Proof.} We follow the proof of \cite{GN}, Proposition 4.11. Put $n={\rm dim}R$. We prove the statement by induction on $n$. Suppose that ${\rm dim}M=n$ and take ${\frak p}\in {\rm Supp}M$ with ${\rm dim}R/{\frak p}=n$. Then ${\rm dim}R_{\frak p}=0$, so that ${\rm depth}R_{\frak p}=0$. Thus ${\frak p}R_{\frak p}\in{\rm Ass}R_{\frak p}$. Hence ${\rm Hom}_{R_{\frak p}}(M_{\frak p},R_{\frak p}/{\frak p}R_{\frak p})\neq 0$ implies ${\rm Hom}_{R_{\frak p}}(M_{\frak p},R_{\frak p})\neq 0$. Thus ${\rm Hom}_{R}(M,R)\neq 0$, i.e., ${\rm grade}M=0$.

When $n=0$, we have ${\rm dim}M=0$. Then the equality holds by above. Let $n>0$. Then we can assume ${\rm dim}M<n$. Since ${\rm dim}R/{\frak p}=n$ for any minimal prime ideal $\frak p$ of $R$, it holds from the assumption that $[0:_RM]\not\subset {\frak p}$ for any minimal prime ideal $\frak p$ of $R$. Thus $[0:_RM]\not\subset{\frak p}$ for any ${\frak p}\in{\rm Ass}R$. Since $[0:_RM]$ is a graded ideal, $[0:_RM]$ contains a homogeneous regular element $x$ by A.1 Proposition. We have that ${\rm Ext}_R^i(M,R)\cong{\rm Ext}_{R/xR}^{i-1}(M,R/xR)$ for $i\geq 0$. Thus ${\rm grade}_{R/xR}M={\rm grade}_RM-1$. By Lemma A.11 and induction, we get ${\rm dim}_{R/xR}M+{\rm grade}_{R/xR}M=n-1$, hence ${\rm dim}_RM+{\rm grade}_RM-1=n-1$, which gives the desired equality. $\square$

\medskip

We state a characterization of a Cohen-Macaulay graded module over a $^*$local ring by means of the $^*$depth and $^*$dimension.

\medskip

A.13. {\sc Theorem}\ \ 
{\sl Let $(R,{\frak m})$ be a $^*$local ring with the condition (P) and $M\in {\rm mod}_0R$. Then $M$ is Cohen-Macaulay if and only if \,$^*{\rm depth}M=\ ^*{\rm dim}M$.}

{\it Proof.}
Put $I=[0:_RM]$ and $\overline{R}=R/I,\ \overline{\frak m}={\frak m}/I$. Then we have that $^*{\rm dim}M={\rm dim}\overline{R}_{\overline{\frak m}}={\rm dim}R_{\frak m}/[0:_{R_{\frak m}}M_{\frak m}]={\rm dim}M_{\frak m}$. It holds from \cite{MR} or \cite{NVO}, Chapter B, Theorem III.2.1 that $M$ is Cohen-Macaulay if and only if $M_{\frak m}$ is Cohen-Macaulay. Look at the following inequalities
$$^*{\rm depth}M={\rm depth}({\frak m},M)\leq {\rm depth}M_{\frak m}\leq{\rm dim}M_{\frak m}=\ ^*{\rm dim}M.$$
If $^*$depth$M=\ ^*{\rm dim}M$, then $M_{\frak m}$ is Cohen-Macaulay by above. Conversely, suppose $M$ to be Cohen-Macaulay. Then depth$({\frak m},M)={\rm depth}M_{\frak m}$ holds by \cite{M}, Theorem 17.3. Thus we get $^*$depth$M=\ ^*{\rm dim}M$ from the above inequalities. $\square$

\medskip

A.14. {\sc Lemma.} (\cite{AB}, Proposition 4.16)\ \ 
{\sl Let $R$ be a commutative Noetherian ring and $X$ a finite $R$-module with \mbox{\rm G-dim}$X<\infty$. Then ${\rm grade}U\geq i$ for all $i>0$ and all $R$-submodules $U$ of ${\rm Ext}_R^i(X,R).$ }

{\it Proof.}
Let ${\frak p}\in{\rm Supp}U$. Then ${\rm Ext}_{R_{\frak p}}^i(X_{\frak p},R_{\frak p})\neq 0$. Hence G-dim$_{R_{\frak p}}X_{\frak p}\geq i$. By Auslander-Bridger formula (\cite{AB}, Theorem 4.13 (b) or \cite{C}, Theorem 1.4.8), it follows that
$${\rm depth}R_{\frak p}={\rm depth}X_{\frak p}+\mbox{G-dim}_{R_{\frak p}}X_{\frak p}\geq \mbox{G-dim}_{R_{\frak p}}X_{\frak p}\geq i.$$
Hence ${\rm grade}U={\rm min}\{{\rm depth}R_{\frak p}:{\frak p}\in{\rm Supp}U\}\geq i$ by \cite{AB}, Corollary 4.6.  $\square$ 

\medskip

A.15. {\sc Lemma.}\ \ 
{\sl Let $R$ be a commutative Noetherian ring and $X$ a finite $R$-module of grade $s$. Assume {\rm G-dim}$X$ to be finite. Then the equality ${\rm grade\,Ext}_R^s(X,R)=s$ holds true.}

{\it Proof.}
When $s=0$, that is, $X^*\neq0$, then $X^{***}\neq0$. Hence $X^{**}\neq0$.

We assume that $s>0$. By A.14, it holds that ${\rm grade\,Ext}_R^s(X,R)\geq s$. The converse inequality follows from \cite{HN}, Lemma 4.4 (Its proof contains trivial misprints : in the last line of p.182, $X_n^*$ should be read $(\Omega^nX)^*$ and three places in line 3-5 of p.183 should be read similarly). Hence we get the desired equality.   $\square$


\begin{thebibliography}{99}


\bibitem{AB}M. Auslander and M. Bridger, {\it Stable module theory}, Mem. of
the  AMS 94, Amer. Math. Soc., Providence 1969.

\bibitem{ARS}M. Auslander, I. Reiten, and S. O. Smalo,{\it Representation Theory of Artin Algebras}, Cambridge stud. Adv. Math. 36, Cambridge Univ. Press, 1995.
\bibitem{B}H. Bass, {\it On the ubiquity of Gorenstein rings}, Math. Z. 82 (1963) 8-28.

\bibitem{Bj1}J-E. Bj\"ork, {\it Rings of Differential Operators}, Math. Ligrary 21, North Holland, 1979.

\bibitem{Bj2}J-E. Bj\"ork, {\it Filtered Noetherian rings}, in Noetherian Rings and Their Applications Ed. L. W. Small, AMS Math, Surveys and Monographs, vol.24, 1987, 59-97.

\bibitem{Bj3}J-E. Bj\"ork, {\it The Auslander condition on Noetherian rings}, in S\'em. d'Alg\`ebre P. Dubreil et M.-P. Malliavin, 1987-88 (M.-P. Malliavin, ed.), Lecture Notes in Math. 1404, Springer, 1989, 137-173.

\bibitem{BE}J-E. Bj\"ork and E. K. Ekstr\"om, {\it Filtered Auslander Gorenstein rings}, in Operator algebras, unitary representations, enveloping algebras, and invariant theory, Progr. Math. 92, Birkhauser (1990) 425-448.


\bibitem{BH}W. Bruns and J. Herzog, {\it Cohen-Macaulay Rings}, Cambridge Stud. Adv. Math. 39, Cambridge Univ. Press, 1993.


\bibitem{C}L. W. Christensen, {\it Gorenstein Dimensions}, Lecture Notes in Mathematics 1747, Springer, 2000.


\bibitem{F}C. Faith, {\it Algebra {\rm I\!I} Ring Theory}, Springer, 1976.


\bibitem{GN}S. Goto and K. Nishida, {\it Towards a theory of Bass numbers with application to Gorenstein algebras}, Colloq. Math. 91(2) (2002) 191-253.


\bibitem{GW}S. Goto and K. Watanabe, {\it On graded rings, {\rm I}}, J. Math. Soc. Japan, 30 (1978) 179-213.


\bibitem{HN}M. Hoshino and K. Nishida, {\it A generalization of the Auslander Formula}, Representations of Algebras and Related Topics, Fields Institute Communications vol. 45 (2005) 175-186.


\bibitem{HVO}L. Huishi and F. Van Oystaeyen, {\it Zariskian Filtrations}, K-Monographs in Mathematics, 2, 1996.

\bibitem{Iw}Y. Iwanaga, {\it Duality over Auslander-Gorenstein rings}, Math. Scand. 81(2) (1997) 184-190. 

\bibitem{I}O. Iyama, {\it Symmetry and duality on $n$-Gorenstein rings}, J. Algebra 269 (2003) 528-535.

\bibitem{K}M. Kashiwara, {\it $D$-modules and microlocal calculus}, Translations of Mathematical Monographs, 217, AMS, 2003.

\bibitem{M}H. Matsumura, {\it Commutative Ring Theory}, Cambridge stud. Adv. Math. 8, Cambridge Univ. Press, 1986.

\bibitem{MR}J. Matijevic and P. Roberts, {\it A conjecture of Nagata on graded Cohen-Macaulay rings}, J. Math. Kyoto Univ., 14 (1974) 125-128.

\bibitem{NVO}C. N\v ast\v asescu and F. Van Oystaeyen, {\it Graded Ring Theory}, Math. Library 28, North Holland, 1982.

\bibitem{R}J. J. Rotman, {\it An introduction to homological algebra}, Academic Press, 1979.

\end{thebibliography}
\end{document}